\documentclass[11pt]{article}% version 20231105

\usepackage{a4,amssymb,amsmath,color,graphicx}
\usepackage{trfsigns}
\usepackage{mathrsfs}
\usepackage{bm}
\usepackage[ruled,vlined]{algorithm2e}
\usepackage[english]{babel}

\usepackage{tikz}
\usetikzlibrary{arrows,shapes}
\tikzstyle{vertex}=[circle,fill=white!25,minimum size=20pt,inner sep=0pt]
\tikzstyle{edge} = [draw,thick,-]
\tikzstyle{edgethin} = [draw,-]
\tikzstyle{weight} = [font=\small]
\usetikzlibrary{decorations}
\usetikzlibrary{decorations.pathmorphing} % (solid) snake :)
\usetikzlibrary{decorations.pathreplacing} % show path construction

\newcommand{\N}{ {\mathbb N} }
\newcommand{\R}{ {\mathbb R} }
\newcommand{\levelN}{\ensuremath{\bm{l}}}
\newcommand{\norm}[1]{\left\Vert #1 \right\Vert}
\newcommand{\gO}{\mathcal{O}}
\newcommand{\ds}{\displaystyle}

\DeclareMathOperator{\supp}{supp}

\NewDocumentCommand{\numKnots}{O{l}}{\ensuremath{m_{#1}}}
\NewDocumentCommand{\knot}{O{l}O{j}}{\ensuremath{y_{#2}^{(#1)}}}
\NewDocumentCommand{\knotHier}{O{l}O{j}}{\ensuremath{x_{#2}^{(#1)}}}
\NewDocumentCommand{\KnotSet}{O{l}}{\ensuremath{Y_{#1}}}
\NewDocumentCommand{\KnotSetHier}{O{l}}{\ensuremath{X_{#1}}}
\NewDocumentCommand{\baseFunc}{O{l}O{j}}{\ensuremath{\phi_{#2}^{(#1)}}}
\NewDocumentCommand{\baseFuncHier}{O{l}O{j}}{\ensuremath{\varphi_{#2}^{(#1)}}}
\NewDocumentCommand{\BaseFuncSet}{O{l}O{j}}{\ensuremath{\varPhi_{#2}^{(#1)}}}
\NewDocumentCommand{\Interpol}{O{l}}{\ensuremath{\mathscr{U}^{(#1)}}}
\NewDocumentCommand{\InterpolDelta}{O{l}}{\ensuremath{\Delta^{(#1)}}}
\NewDocumentCommand{\InterpolSmol}{O{w}}{\ensuremath{\mathcal{A}_{#1}}}
\NewDocumentCommand{\InterpolI}{O{w}}{\ensuremath{\mathcal{I}^{(#1)}}}
\NewDocumentCommand{\leveld}{O{d}}{\ensuremath{l_{#1}}}
\NewDocumentCommand{\weight}{O{l}O{j}}{\ensuremath{w_{#2}^{(#1)}}}
\NewDocumentCommand{\support}{O{l}O{j}}{\ensuremath{s_{#2}^{(#1)}}}
\NewDocumentCommand{\Stencil}{O{x}}{\ensuremath{\mathcal{S}_{#1}}}
\NewDocumentCommand{\fcurve}{O{N}}{\ensuremath{f^{#1}_{\text{0}}}}
\NewDocumentCommand{\fc}{O{N}}{\ensuremath{f^{#1}_{\text{1}}}}
\NewDocumentCommand{\fg}{O{N}}{\ensuremath{f^{#1}_{\text{2}}}}
\NewDocumentCommand{\fd}{O{N}}{\ensuremath{f^{#1}_{\text{3}}}}
\NewDocumentCommand{\errMax}{}{\ensuremath{\epsilon_{\infty}}}
\NewDocumentCommand{\errTwo}{}{\ensuremath{\epsilon_{2}}}

\newtheorem{remark}{Remark}[section]

\begin{document}
\title{Adaptive hp-Polynomial Based \\
Sparse Grid Collocation Algorithms for \\
Piecewise Smooth Functions with Kinks}
\author{Hendrik Wilka and Jens Lang \\
{\small \it Technical University Darmstadt,
Department of Mathematics} \\
{\small \it Dolivostra{\ss}e 15, 64293 Darmstadt, Germany} \\
{\small wilka@mathematik.tu-darmstadt.de, lang@mathematik.tu-darmstadt.de}}
\maketitle

\begin{abstract}
High-dimensional interpolation problems appear in various applications
of uncertainty quantification, stochastic optimization and machine learning.
Such problems are computationally expensive and request the use of adaptive 
grid generation strategies like anisotropic sparse grids to mitigate the curse of 
dimensionality. However,
it is well known that the standard dimension-adaptive sparse grid method 
converges very slowly or even fails in the case of non-smooth functions. 
For piecewise smooth functions with kinks, we
construct two novel $hp$-adaptive sparse grid collocation algorithms that combine
low-order basis functions with local support in parts of the domain 
with less regularity and variable-order basis functions elsewhere. Spatial refinement
is realized by means of a hierarchical multivariate knot tree which allows the
construction of localised hierarchical basis functions with varying order. 
Hierarchical surplus is used as an error indicator to automatically
detect the non-smooth region and adaptively refine the collocation points there. 
The local polynomial degrees are optionally selected by a greedy approach or a kink
detection procedure. Four numerical benchmark examples with different dimensions are 
discussed and comparison with locally linear, quadratic and highest degree basis 
functions are given to show the efficiency and accuracy of the proposed methods.
\end{abstract}


\section{Introduction}
The efficient approximation of high-dimensional functions has become a fundamental component in
almost all scientific fields. Prominent areas are uncertainty quantification, stochastic optimization 
and machine learning. For example, the quantitative characterization of uncertainties in certain
outcomes, if some input parameters are stochastically modelled, typically leads to high-dimensional
problem spaces in which interpolation operations need to be performed. The use of full tensor grids
is often not feasible since the number of grid points increases exponentially, i.e., $O(M^N)$ for
$M$ grid points in each of the $N$ dimensions -- the well-known curse of dimensionality. Sparse grids
have emerged as an extremely useful tool box offering the possibility to only select grid points that
contribute most to the multi-dimensional approximation of the solution 
\cite{GerstnerGriebel1998,Smolyak1963,Zenger1991}. The complexity of the sparse grid method
reduces to $O(M\log(M)^{N-1})$ \cite{BungartzGriebel2004}. By taking advantage of higher
smoothness and dimension-dependent adaptivity, which was originally proposed in \cite{GerstnerGriebel2003} 
as anisotropic grids and further developed in \cite{Klimke2006,NobileTemponeWebster2008}, even 
faster rates of convergence can be achieved. Given appropriate error estimators, this generalized 
sparse grid method automatically concentrates grid points in important dimensions.

However, when the function to be approximated exhibits a non-smooth dependence on the uncertain 
input parameters, sparse grid methods that use a global polynomial basis can not efficiently resolve 
local kinks or discontinuities in the random space. In such cases, spatial refinement has to be added
to boost the efficiency of sparse grid approximations. Locally adaptive sparse grids were first 
investigated in \cite{Griebel1998} for the numerical solution of PDEs. So-called $h$-adaptive 
generalized sparse grid methods have been used in \cite{MaZabaras2009} with multi-linear hierarchical 
basis functions of local support and in \cite{JakemanRoberts2011} with a higher, but fixed maximum 
polynomial degree.
The magnitude of the hierarchical surplus, i.e. the difference between the function and its 
approximation at a newly selected grid point not yet part of the current point set, guides the 
spatial refinement. More
recently, a generalized spatially adaptive sparse grid combination technique with dimension-wise refinement
was proposed in \cite{ObersteinerBungartz2021}. Simplex stochastic collocation for piecewise 
smooth functions with kinks was analysed in \cite{FuchsGarcke2020}.
Although, all these methods achieve a fast rate of convergence in smooth regions and acceptable 
accuracies around discontinuities or kinks, a combination of local $h$- and $p$-adaptivity in 
the spirit of the $hp$-version of the finite element method which often shows exponential 
convergence has not been utilised so far. 
 
Thus, we exploit the localised polynomial basis proposed in \cite{Bungartz1998} to construct an 
$hp$-adaptive generalized sparse grid method. Local polynomial basis functions of different 
degree can be constructed by means of so-called ancestor knots. With the help of the hierarchical 
surplus, we adopt the local polynomial degree to the local regularity of the function. We have 
developed two approaches: (1) The greedy approach derives an individual score for all possible 
polynomial degrees and chooses the one with the minimal
interpolation error calculated at newly selected collocation points. (2) The kink detection approach based on
a polynomial annihilation technique \cite{ArchibaldGelbYoon2005,ArchibaldGelbYoon2008} is applied to identify
regions with discontinuous first derivatives. Then, piecewise linear basis functions are used there. 
New collocation points are investigated locally around an active set of already chosen points and added to that set
until a given threshold for the hierarchical surplus is reached.   

The paper is organised as follows. In Section~\ref{sec:ahsg}, we introduce the fundamentals of 
adaptive hierarchical sparse grids based on localised hierarchical basis functions with varying order.
Our novel $hp$-generalised sparse grid algorithm is described in Section~\ref{sec:hpgsg}. Numerical results
for four benchmark functions and comparison to other approaches are presented in Section~\ref{sec:num}. We
summarize our results in Section~\ref{sec:concl}.

\section{Adaptive hierarchical sparse grids}\label{sec:ahsg}
We are interested in interpolating non-smooth functions $f:\Omega\rightarrow\R$ defined on an $N$-dimensional
bounded domain $\Omega$. Without loss of generality, we choose $\Omega=[-1,1]^N$ and require that $f$ can be evaluated at
arbitrary points in $\Omega$. Our main focus will be on piecewise smooth functions that have kinks along certain manifolds.

\subsection{From univariate interpolation to hierarchical sparse grids} 
Let us start with univariate interpolation on the interval $I = [-1,1]$.
For a given index $l\in\N_0$, we consider $\numKnots$ knots $-1 \le \knot[l][1] < ... < \knot[l][\numKnots] \le 1$ with 
the corresponding set of knots 
\begin{align*}
\KnotSet[l] := \left\{\knot[l][1], ..., \knot[l][\numKnots] \right\}, \quad l \in \N_0.
\end{align*}
We assume $\numKnots[l+1] \ge \numKnots$ for different levels of interpolation. A knot-based interpolation operator $\Interpol$ 
is now given by
\begin{align*}
\Interpol(f)(x) = \sum_{j=1}^{\numKnots} f(\knot) \baseFunc(x), \quad x \in I,
\end{align*}
where $\baseFunc$, $j=1,...,\numKnots$, are basis functions which satisfy $\baseFunc(\knot[l][k]) = \delta_{jk}$. Hence, 
$\Interpol(f)$ is exact on the given interpolation knots.

Defining the multiindex $\levelN = (\leveld[1],...,\leveld[N])$ and 
$\KnotSet[\levelN] = \KnotSet[\leveld[1]] \times ... \times \KnotSet[\leveld[N]]$ with knot sets 
$\KnotSet[\leveld[1]],...,\KnotSet[\leveld[N]]$ for each dimension, we get the multivariate interpolation operator
\begin{align}
\Interpol[\levelN](f)(\bm{x}) &= \left( \Interpol[\leveld[1]] \times ... \times \Interpol[\leveld[N]] 
\right)(f)(\bm{x}) \nonumber \\
&= \sum_{j_1 = 1}^{\numKnots[\leveld[1]]} ... \sum_{j_N = 1}^{\numKnots[\leveld[N]]}
f\left( \knot[\leveld[1]][j_1] ,...,\knot[\leveld[N]][j_N] \right) \left( \baseFunc[\leveld[1]][j_1] \cdot ... \cdot \baseFunc[\leveld[N]][j_N] \right)(\bm{x}), \quad \bm{x}\in\Omega. \label{fullTensorGrid} 
\end{align}
Again, this construction satisfies the interpolation property $\Interpol[\levelN](f)(\knot[\levelN][\bm{j}]) = f(\knot[\levelN][\bm{j}])$ for all $\knot[\levelN][\bm{j}]=(\knot[\leveld[1]][j_1] ,...,\knot[\leveld[N]][j_N]) \in \KnotSet[\levelN]$.
Since the number of interpolation knots in the full-tensor product \eqref{fullTensorGrid} grows very quickly, we apply a sparse grid interpolation
based on the Smolyak algorithm \cite{Smolyak1963}. Assuming nested point sets $\KnotSet \subset \KnotSet[l+1]$ for $l \in \N_0$, we introduce new
point sets
\begin{align*}
\KnotSetHier[0] &:= \KnotSet[0],\quad \KnotSetHier[l] := \KnotSet[l] \backslash \KnotSet[l-1], \;l \ge 1.
\end{align*}
With $\KnotSetHier = \left\{ \knotHier[l][1], ... , \knotHier[l][\left| \KnotSetHier \right|] \right\}$, we have $\knotHier = \knot[l][j']$ for some $j'$. On the other hand, for every $\knot$ there exist a unique pair $(l',j')$ such that $\knot = \knotHier[l'][j']$, 
since $\KnotSet = \cup_{k=1,...,l} \KnotSetHier[k]$. Further, the corresponding incremental interpolation operators 
\begin{align*}
\InterpolDelta[l] = \Interpol[l] - \Interpol[l-1],\;l\in\N_0,
\end{align*}
with the convention $\Interpol[-1] \equiv 0$ has the property
\begin{align*}
\InterpolDelta(f)(\knotHier) &= \Interpol(f)(\knotHier) - \Interpol[l-1](f)(\knotHier)
= f(\knotHier) - \Interpol[l-1](f)(\knotHier)
\end{align*}
and therefore $\InterpolDelta(f)(\knotHier[l'][j'])=0$ for $l'<l$. Following \cite[Sec. 3.3]{MaZabaras2009}, we get the well-known 
relation
\begin{align}
\InterpolDelta(f) &=  \sum_{\knotHier \in \KnotSetHier} \left(  f(\knotHier) - \Interpol[l-1](f)(\knotHier)\right) 
\baseFuncHier =: \sum_{\knotHier \in \KnotSetHier} \weight \baseFuncHier. \label{DifferenceFormula}
\end{align}
Here, $\weight$ is the one-dimensional hierarchical surplus representing the difference between the function values at the current
and the previous interpolation levels. The set of functions $\baseFuncHier$ are now defined as the hierarchical basis functions.

Using the identity $\sum_{k=0,\dots,l} \InterpolDelta = \Interpol[l]$, the full-tensor approximation \eqref{fullTensorGrid} can be
decomposed into
\begin{align}
\Interpol[\levelN] &= \Interpol[\leveld[1]] \times ... \times \Interpol[\leveld[N]] =
\left( \sum_{\leveld[1] = 0}^{w} \InterpolDelta[\leveld[1]]  \right)
\times ... \times 
\left( \sum_{\leveld[N] = 0}^{w} \InterpolDelta[\leveld[N]]  \right) \nonumber \\
&=\sum_{\leveld[1] = 0}^{w} ... \sum_{\leveld[N] = 0}^{w} \InterpolDelta[\leveld[1]] 
\times ... \times 
\InterpolDelta[\leveld[N]] =: \sum_{\levelN:\,\norm{\levelN}_{\infty} \le w} \InterpolDelta[\levelN]
\end{align}
for any $w\in\N_0$. Defining $L:=|\bm{l}|=\sum_{d=1,\ldots,N}l_d$, the classical sparse grid interpolation operator
originally introduced in \cite{Smolyak1963} is given by 
\begin{align*}
\InterpolSmol[q](f) =\sum_{\levelN:\, L \le q} \InterpolDelta[\leveld[1]](f) 
\times ... \times \InterpolDelta[\leveld[N]](f)
\end{align*}
for some maximum order $q\ge N$.
Equivalently, we can use Smolyak's formula \cite{GerstnerGriebel1998,Smolyak1963}
\begin{align}
\InterpolSmol[q](f) =  \sum_{\levelN: q - N + 1 \le L \le q} (-1)^{q - L} 
\begin{pmatrix}
N - 1 \\ q - L
\end{pmatrix}
\left( \Interpol[\leveld[1]](f) \times ... \times \Interpol[\leveld[N]](f) \right) \label{SmolyakScheme}
\end{align}
to get a representation without incremental operators. For $M$ knots per dimension, we have $M^N$ function
evaluations in the full-tensor approximation \eqref{fullTensorGrid}, which grows exponentially with the
dimension $N$ and hence suffers from the curse of dimensionality. Sparse grids reduce the number of points
to $\gO(M \log(M)^{N-1})$ while keeping a comparable approximation quality for smooth functions \cite{BungartzGriebel2004}.
Obviously, the curse of dimensionality is significantly delayed, but not broken. Further improvements can be made by adaptive
strategies, i.e., only adding knots that reduce the interpolation error efficiently \cite{GerstnerGriebel2003,NobileTemponeWebster2008}.
However, problems with non-smooth functions $f$ require spatial refinement and the use of low order basis functions
in parts of the domain with less regularity. A generalized sparse grid approach which allows such an improvements will be described next.

\subsection{Hierarchical basis functions with local support}
In order to identify and resolve local non-smooth variations of the function $f$, detection algorithm and spatial refinement 
can be performed on the level of a single hierarchical basis function $\baseFuncHier$ from \eqref{DifferenceFormula}. We consider
one-dimensional equidistant points of the sparse grid which can be considered as a tree-like data structure 
\cite{Bungartz1998,BungartzGriebel2004}. Let an equidistant knot set $Y_l$ on $[-1,1]$ given by
\begin{alignat}{3}
\knotHier[0][1] &= 0, && \nonumber \\
\knotHier[1][j] &= 2j-3,\quad &&\,j=1,2=: |X_1| \nonumber\\
\knotHier &=  (2j -1) 2^{1-l} - 1,\quad &&\,j=1,...,2^{l-1} =: |X_l|,\,l\ge 2. \label{hier:nodeset}
\end{alignat}
For these knots, we can define a hierarchical structure: We start with a central knot, $\knotHier[0][1]\!=\!0$. This knot has the two child knots $\knotHier[1][1]\!=\!-1$ and $\knotHier[1][2]\!=\!1$ located at the boundary of the interval. These knots have the child knots $\knotHier[2][1]\!=\! -\frac{1}{2}$ and $\knotHier[2][2]\!=\!\frac{1}{2}$, respectively. From now on, the child knots can be defined recursively. Every knot $\knotHier[l][j]$ with $l \ge 2$ has two child knots with the positions $\knotHier[l][j] \pm 2^{- l}$. In the following, $S_C(\knotHier)$ denotes the children set of a knot $\knotHier$. Further, let $S_P(\knotHier)$ be the parent of a knot, which always exists for $l\ge 1$. We have
\begin{align*}
\knotHier[l-1][j'] = S_P(\knotHier) \Leftrightarrow \knotHier \in S_C(\knotHier[l-1][j']).
\end{align*}
With that, we can also define the set of unique ancestors $S_A(\knotHier):=\{\knotHier[l-1][j_{l-1}],\knotHier[l-2][j_{l-2}],\ldots,\knotHier[0][1]\}$ that satisfy
\begin{align*}
S_P(\knotHier) = \knotHier[l-1][j_{l-1}],\, S_P(\knotHier[l-1][j_{l-1}]) = \knotHier[l-2][j_{l-2}],\ldots , 
S_P(\knotHier[1][j_{1}]) = \knotHier[0][1].
\end{align*}
The resulting structure forms a knot tree, see Fig.~\ref{fig:hier_tree} for an illustration with $l=4$.
\begin{figure}[h]
	\begin{center}
		\begin{tikzpicture}
			\node[vertex](1) at (4,4){$0$};
			\node[vertex](2) at (0,3){$-1$};
			\node[vertex](3) at (8,3){$1$};
			\node[vertex](4) at (2,2){$-\frac{1}{2}$};
			\node[vertex](5) at (6,2){$\frac{1}{2}$};
			\node[vertex](6) at (1,1){$-\frac{3}{4}$};
			\node[vertex](7) at (3,1){$-\frac{1}{4}$};
			\node[vertex](8) at (5,1){$\frac{1}{4}$};
			\node[vertex](9) at (7,1){$\frac{3}{4}$};
			\node[vertex](10) at (0.5,0){$-\frac{7}{8}$};
			\node[vertex](11) at (1.5,0){$-\frac{5}{8}$};
			\node[vertex](12) at (2.5,0){$-\frac{3}{8}$};
			\node[vertex](13) at (3.5,0){$-\frac{1}{8}$};
			\node[vertex](14) at (4.5,0){$\frac{1}{8}$};
			\node[vertex](15) at (5.5,0){$\frac{3}{8}$};
			\node[vertex](16) at (6.5,0){$\frac{5}{8}$};
			\node[vertex](17) at (7.5,0){$\frac{7}{8}$};
			% Connect
			\path[edgethin] (1) -- (2);
			\path[edgethin] (1) -- (3);
			\path[edgethin] (2) -- (4);
			\path[edgethin] (3) -- (5);
			\path[edgethin] (4) -- (6);
			\path[edgethin] (4) -- (7);
			\path[edgethin] (5) -- (8);
			\path[edgethin] (5) -- (9);
			\path[edgethin] (6) -- (10);
			\path[edgethin] (6) -- (11);
			\path[edgethin] (7) -- (12);
			\path[edgethin] (7) -- (13);
			\path[edgethin] (8) -- (14);
			\path[edgethin] (8) -- (15);
			\path[edgethin] (9) -- (16);
			\path[edgethin] (9) -- (17);
		\end{tikzpicture}
    \parbox{14cm}{
	\caption{The first $5$ levels of the equidistant knot tree on $[-1,1]$. For example,
the ancestors of knot $x=-\frac58$ are $S_A(-\frac58)=\{-\frac34,-\frac12,-1,0\}$.}
	\label{fig:hier_tree}
    }
    \end{center}
\end{figure}

We will now define designated supports for local hierarchical basis functions. Let
$\support\!:=\!\supp\baseFuncHier$ and set
\begin{alignat}{2}
\support[0][1] &= [-1,1],\; \support[1][1] = [-1,0],\; \support[1][2] = [0,1], \nonumber\\[1mm]
\support[l][j] &= \left[ \knotHier - 2^{1-l}, \knotHier + 2^{1-l} \right],\;l\ge 2. 
\label{suppHier}
\end{alignat}
This yields
\begin{align*}
\support[l][j_1] \cap \support[l][j_2] = 
\begin{cases}
\left\{ \knotHier[l'][j']\right\} \text{ with } l' < l ,& \text{if } |j_1 - j_2| = 1 \\[2mm]
\emptyset ,& \text{if } |j_1 - j_2| > 1 \\
\end{cases}.
\end{align*}
Thus, supports of knots in the same level intersect at most in one point. Additionally, the supports 
are nested and respect the hierarchical structure of the knots, i.e.,
\begin{align*}
\knotHier \in S_C(\knotHier[l-1][j']) \Rightarrow \support \subset \support[l-1][j'].
\end{align*}
A new basis function $\baseFuncHier[l][j]$ modifies an existing interpolation only in the support of its parent and 
$\baseFuncHier[l][j]$, $\baseFuncHier[l][j']$ with $j\ne j'$ do not influence each other.

The parent-child relation for knots can be extended to the multivariate case by defining the sets of children
\begin{align*}
S_C(\knotHier[\levelN][\bm{j}]) := \left\{  \left( \knotHier[\leveld[1]][j_1] ,\ldots , \knotHier[\leveld[d-1]][j_{d-1}], x,  
\knotHier[\leveld[d+1]][j_{d+1}],\ldots ,\knotHier[\leveld[N]][j_{N}] \right) : x \in S_C(\knotHier[l_d][j_d]) 
\text{ for } d=1,\ldots,N  \right\}
\end{align*}
for $\levelN=(l_1,\ldots,l_N)$ and $\bm{j}=(j_1,\ldots,j_N)$. Observe that each multivariate knot can have at most $2N$ children
and multiple parents along different axes.

Local polynomial basis functions of high degree can be constructed by means of ancestor knots as originally studied in
\cite{Bungartz1998}. Let $N\!=\!1$ and $l>\!2$. Then, we use knots from $S_A(\knotHier)$ to define a higher-order basis 
function $\baseFuncHier[l][j]$ of the difference operator $\InterpolDelta$ associated to the knot $\knotHier[l][j]$ and 
satisfying 
\begin{align*}
\baseFuncHier[l][j](\knotHier[l][j])=1,\quad\baseFuncHier[l][j](\knotHier[l'][j']) &= 0 \text{ for all } 
\knotHier[l'][j'] \in S_A(\knotHier).
\end{align*}
We note that the boundary points of the support $\support$ of $\baseFuncHier[l][j]$ are also included in $S_A(\knotHier)$.
Further ancestors lie outside of $\support$. The definition of the local higher-order basis functions is consistent in the sense
that a global interpolant built up with them fulfills all requirements of a hierarchical Lagrangian interpolation. The
new weight is given by $\weight=f(\knotHier)-\InterpolI[l-1](\knotHier)$, where $\InterpolI[l-1]$ denotes the interpolation 
operator that uses knots up to level $l-1$. The number of ancestors knots used defines the degree of $\baseFuncHier[l][j]$.
To reach order $p$, we take the knot $\knotHier$, its two ancestors $a_1,a_2$ at the boundary of $\support$, and additionally
its nearest $p-2$ ancestors $a_3,\ldots,a_p$ ordered by their distance to $\knotHier$. The new basis function is then defined
by
\begin{align}
\baseFuncHier[l][j](x):=\baseFuncHier[l][j,p](x) = 
\begin{cases}
\ds\prod_{k=1}^{p} \frac{x - a_k}{\knotHier - a_k} &,\, x \in \support, \\[2mm]
0 &,\, x \notin \support.
\end{cases} 
\label{baseFunctionP}
\end{align}
With a predefined maximum degree $p_{\max}$, we will consider the set of basis functions
\begin{align}
\BaseFuncSet &:= \left\{ \baseFuncHier[l][j,p],\, p=2,\ldots,\min\{p_{\max}, l\} \right\}.
\label{setOfBaseFunc}
\end{align}
It remains to choose basis functions for $l\le 2$. For the root knot $\knotHier[0][1]$, we
take the constant function $\baseFuncHier[0][1,0] \equiv 1$. In the case $l\!=\!1$, linear basis
functions $\baseFuncHier[1][1,1](x)$, $\baseFuncHier[1][2,1](x)$ on $\support[1][1]$,
$\support[1][2]$ defined in \eqref{suppHier} are chosen, respectively. Eventually, the quadratic
polynomial $\baseFuncHier[2][j,2](x)$ vanishing at the boundary of the support is taken for $l\!=\!2$. But, we will
also use the piecewise linear basis function $\baseFuncHier[2][j,1](x)$ as defined for $l\!=\!1$.

Generalized multivariate basis functions of degree $\bm{p}=(p_1,\ldots,p_N)\in\N^N$ are now defined as
products 
\begin{align}
\baseFuncHier[\bm{l}][\bm{j},\bm{p}](\bm{x}) &:= \baseFuncHier[l_1][j_1,p_1](x_1)\cdot\ldots\cdot
\baseFuncHier[l_N][j_N,p_N](x_N)
\label{multBasisFunc}
\end{align}
of one-dimensional basis functions, where $\bm{p}$ satisfies $p_i\le\min\{p_{\max}, l_i\}, i=1,\ldots,N$, 
and $|\bm{l}|\le q\in\N$. We note that resulting global interpolants are continuous functions on $\Omega$.
\begin{remark}
In \cite{Bungartz1998}, the use of one $p:=p_{\max}$ for all directions is called the p-regular scenario.
For functions $f$ having bounded mixed derivatives up to order $p+1$ in the $L^q$-norm, $q\in\{2,\infty\}$, 
and vanishing on 
the boundary of $\bar\Omega:=[0,1]^N$, the $L^\infty$- and $L^2$-error of the p-regular sparse grid 
approximation is of order $\gO(h_n^{p+1}|\log_2h_n|^{N-1})$, $h_n=2^{-n}$, if $n\ge p$ 
\cite[Theorem 4.8]{BungartzGriebel2004}. This shows a higher-order approximation.
\end{remark}
Multi-linear hierarchical basis functions, i.e. $p=1$ in all directions, of local support have been successfully applied
in \cite{MaZabaras2009}. Especially for discontinuous functions, this approach leads to a significant
reduction in the number of collocation points required to achieve comparable accuracy as classical,
possibly dimension-adaptive sparse grid collocation methods. Linear adaptive methods, however, put
often an unnecessarily large number of points in region where the solution is smooth. An $h$-adaptive
generalised sparse grid method ($h$-GSG) with fixed maximum degree $p_{\max}$, which is based on the above localised
polynomial basis and the magnitude of the hierarchical surplus, has been proposed in \cite{JakemanRoberts2011}.
It achieves a fast rate of convergence in smooth regions and accuracies around discontinuities or kinks that
are comparable to those of linear methods.
\begin{remark}
An error analysis in \cite{JakemanRoberts2011} reveals that if 
$\|f-\mathscr{U}_{opt}(f)\|\le \varepsilon$ with the least number of function evaluations to
construct the (optimal) interpolant $\mathscr{U}_{opt}(f)$, then the $h$-adaptive GSG approximation
$\mathscr{U}_{N}(f)$ satisfies $\|f-\mathscr{U}_{N}(f)\|\leq \varepsilon(1+C(\varepsilon))$, where
$C(\varepsilon)$ is the number of points that are used in $\mathscr{U}_{opt}(f)$, but not in $\mathscr{U}_{N}(f)$.
$C(\varepsilon)$ depends on the smoothness of $f$: the smoother $f$ the smaller $C(\varepsilon)$ will be.
\end{remark}
\subsection{Localised hierarchical basis functions with varying order}
Although the support $\support[l][j]$ of the one-dimensional basis function $\baseFuncHier[l][j,p]$ gets smaller for
larger $l$ due to its definition in \eqref{suppHier}, ancestors of $\knotHier[l][j]$ also lie outside for $p>2$. This
bears difficulties in the approximation quality when higher-order Lagrangian interpolation is applied across 
regions of less regularity as discontinuities or kinks of the function to be approximated. 

\begin{figure}[t]
\centering
\includegraphics[width=11cm]{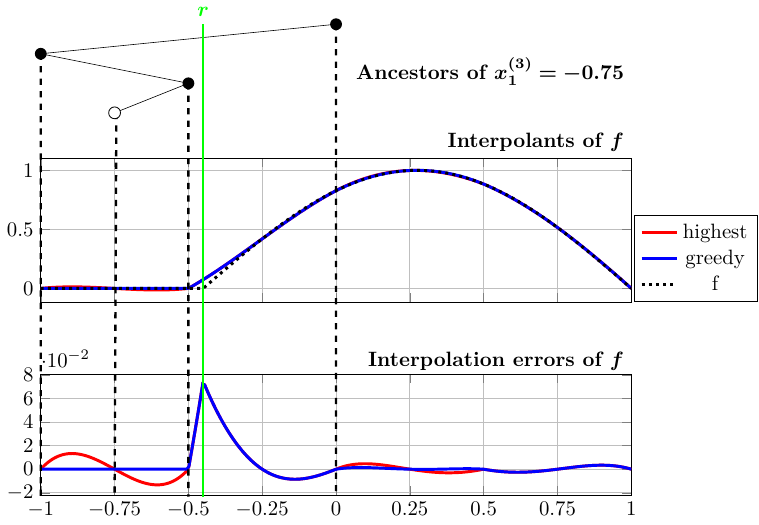}
\parbox{14cm}{
\caption{Interpolation across a kink for $f$ defined in \eqref{funcWithKink}. The basis function
$\baseFuncHier[3][1,3]$ at knot $\knotHier[3][1]=-0.75$ with \textit{highest} possible polynomial degree 
$p=3$ uses a function value at $x=0$, which leads to a relatively bad approximation of $f(x)\equiv 0$ in 
its support $[-1,-0.5]$ away from the kink at $x=-0.45$. The \textit{greedy} approach described in 
Section~\ref{sec:greedy} restricts the used ancestors and yields a zero error there.}
\label{fig:kink1d}
}
\end{figure}

For an illustration, let us consider the function
\begin{align}
f(x) =&\,\left\{ 
\begin{array}{ll}
0 & ,\,x\le r\\
\ds\sin\left( \frac{x-r}{1-r}\pi \right) & ,\, x>r
\end{array}\right.
\label{funcWithKink}
\end{align}
with a kink at $r=-0.45$. From Fig.~\ref{fig:kink1d} it becomes clear that the basis function
for the knot $\knotHier[3][1]=-0.75$ at level $l=3$ uses function values from ancestor knots lying on both sides 
of the kink, which yields an unsatisfactory approximation in its support $[-1,-0.5]$. Adding further
points near the kink will only slowly reduce the interpolation error. Detecting the kink at
an early stage by e.g. a greedy approach explained in Section~\ref{sec:greedy} overcomes this
undesired behaviour. Here, basis functions with support away from the kink are constructed
using knots on one side of the kink only.

\section{The $hp$-generalised sparse grid algorithm}\label{sec:hpgsg}
Here, we will departure 
from previous work in the literature and use localised hierarchical basis functions with varying order. It
is our goal to combine the strengths of both local $h$- and $p$-adaptivity in a generalised
sparse grid method shortly denoted by $hp$-GSG.

Let $A$ denote the active set of adaptively chosen knots $\knotHier[\bm{l'}][\bm{j'}]$ used at level 
$q$, $q\ge |\bm{l'}|$ for all $\bm{l'}$ selected, to construct the interpolant $\mathscr{U}_{A}^{(q)}(f)$.
Then the weight at each child knot $\knotHier[\bm{l}][\bm{j}]\in S_C(\knotHier[\bm{l'}][\bm{j'}])$ 
with $|\bm{l'}|=q$ can be computed by
\begin{align*}
\weight[\bm{l}][\bm{j}] =&\, f(\knotHier[\bm{l}][\bm{j}]) - \mathscr{U}_{A}^{(q)}(f)(\knotHier[\bm{l}][\bm{j}]).
\end{align*}
The new function values $f(\knotHier[\bm{l}][\bm{j}])$ can be used to adopt the current polynomial degree of the 
previously defined basis functions $\baseFuncHier[\bm{l'}][{\bm{j'}}]$, $|\bm{l'}|\!=\!q$. This is done in the
subroutine \textit{Modification}. If a child knot is added to the active set $A$, i.e., we find
$|\weight[\bm{l}][\bm{j}]| \ge w_{\max}$ with a certain threshold $w_{\max}$, a new basis function 
$\baseFuncHier[\bm{l}][{\bm{j}}]$ with an
appropriate polynomial degree is defined in the subroutine \textit{Selection}. Both routines will be explained
in the next subsections. The whole structure of the novel $hp$-GSG algorithm is outlined in Algorithm~\ref{AlgoSparseGridMain}. For
a better readability, we use \glqq evaluate\grqq{} for a function evaluation at a certain knot.
\begin{algorithm}[t]
	\SetAlgoLined
	\KwIn{Function $f(x)$, minimum level $q_{\min}$, maximum level $q_{\max}$, threshold $w_{\max}$, maximum polynomial degree $p_{\max}$}
	\KwResult{Grid $\mathcal{G}$ with knots $\knotHier[\bm{l}][\bm{j}]$, weights $\weight[\bm{l}][\bm{j}]$, basis functions $\baseFuncHier[\bm{l}][\bm{j}]$}
	\Begin{
		Set $\knotHier[\bm{0}][\bm{1}]:=\bm{0}$ and evaluate it.\\
		Set  $A_0:=\{\knotHier[\bm{0}][\bm{1}]\}$, $\mathcal{G}=A_0$, $q:=1$.\\[1mm]
		\While{$q \le q_{\max}$}{
			Set the child set $C = \emptyset$.\\[1mm]
			\For{Every knot $\knotHier[\bm{l'}][\bm{j'}] \in A_{q-1}$}{
				Add every child knot in $S_C(\knotHier[\bm{l'}][\bm{j'}])$ to $C$ and evaluate them.\\
				Define $\baseFuncHier[\bm{l'}][{\bm{j'}}]:=$ \sc{Modification}($\bm{l'},\bm{j'}$).\\			
			}
			\For{Every child $\knotHier[\bm{l}][\bm{j}] \in C$}{
				Calculate the weight $\weight[\bm{l}][\bm{j}]$ for $\knotHier[\bm{l}][\bm{j}]$.\\
				\If{$q \le q_{\min}$ or $|\weight[\bm{l}][\bm{j}]| \ge w_{\max}$}{
				Add $\knotHier[\bm{l}][\bm{j}]$ to $A_q$.\\
				Define $\baseFuncHier[\bm{l}][{\bm{j}}]:=$ \sc{Selection}($\bm{l},\bm{j}$).	
				}	
			}
			Set $\mathcal{G} := \mathcal{G} \cup A_q$\\
		}
	}
	\caption{$hp$-generalised sparse grid approximation}
	\label{AlgoSparseGridMain}
\end{algorithm}

\subsection{Greedy approach}\label{sec:greedy}
For a knot $\knotHier[\bm{l}][\bm{j}]=(\knotHier[l_1][j_1],\ldots,\knotHier[l_N][j_N])\in A$,
the one-dimensional basis function  $\baseFuncHier[l_d][j_d,p_d](x_d)$ used in \eqref{multBasisFunc}
has to be chosen from the set $\BaseFuncSet[l_d][j_d]\cup\{\baseFuncHier[l_d][j_d,1]\}$ defined 
in \eqref{setOfBaseFunc} and thereafter. We
want to optimize all previously chosen $p_d$ in such a way that the interpolation error at some test knots is as small
as possible. Whenever new child knots are considered in the $hp$-GSG method, the two new function
values for varying $\knotHier[l_d][j_d]$, which have to be calculated in any case, can be used to derive a
score for all possible polynomial degrees. Remembering $\baseFuncHier[l_i][j_i,p_i](\knotHier[l_i][j_i])=1$,
$1\le i\le N$, we set
\begin{align*}
S_p :=&\, \max_{x \in S_C(\knotHier[l_d][j_d])} \left| \Interpol[l_d-1] (x) + \weight[{l_d}][j_d] 
\baseFuncHier[l_d][j_d,p](x) - f(x)  \right|,\quad p=1,\ldots,\max\{p_{\max},l_d\},
\end{align*} 
and choose a new $p_d:=\arg \min_{p} S_p$. Note that $\Interpol[l_d-1] + \weight[{l_d}][j_d]\baseFuncHier[l_d][j_d,p_d]$ 
is the local interpolation in level $l_d$ if one would choose $p_d$ as polynomial degree. 
Finally, $N$ individual optimization steps are performed to find $\bm{p}=(p_1,\ldots,p_N)$. This greedy-like strategy
is implemented in \textit{Modification} in Algorithm~\ref{AlgoSparseGridMain}. 

In the subroutine \textit{Selection}, appropriate polynomial degrees have to be chosen for basis functions $\baseFuncHier[\bm{l}][\bm{j}]$
defined for a selected child knot $\knotHier[\bm{l}][\bm{j}]\!\in\! A_q$ on the highest level $q\!=\!|\bm{l}|$. Suppose 
$\knotHier[\bm{l}][\bm{j}]$ was added by changing $\knotHier[l'_d][j'_d]$ in its parent knot $\knotHier[\bm{l'}][\bm{j'}]$
along the $d$-axis. Then we set $p_d:=\min\{p'_d+1,p_{\max}\}$ and replace $\baseFuncHier[l'_d][j'_d,p'_d](x_d)$ by
$\baseFuncHier[l_d][j_d,p_d](x_d)$ in $\baseFuncHier[\bm{l'}][\bm{j'}]$ to define the new basis function for the child knot.

With these definitions of the subroutines, the $hp$-GSG is uniquely determined. For later use, we call it $hp$-GSG-g, where
\glqq g\grqq{} stands for \textbf{g}reedy.

\subsection{Kink detection approach}
Here, we would like to combine our $hp$-GSG method with a kink detection procedure as described in \cite{ArchibaldGelbYoon2005,ArchibaldGelbYoon2008}. Although the greedy approach is in principle able to detect any
unsmoothness of the function $f$, a kink detection might be preferable if it is a priori known that $f$ exhibits kinks, i.e.,
discontinuities in its first derivatives.

We will first recapitulate the main points of the derivative discontinuity detection method \cite{ArchibaldGelbYoon2005} and
then describe its incorporation into our sparse grid algorithm. It is sufficient to consider the one-dimensional case and
apply the procedure on each axis separately. 

Let $f$ be a continuous function with well defined left and right side limits. Assume that $f \in C^1(I \setminus \Theta)$
with $I:=[-1,1]$ and $\Theta$ being the discrete set of jump discontinuities in the first derivative. For each $x \in I$, 
we define the usual jump function
\begin{align*}
\left[ f' \right](x) := \lim_{x_+ \searrow x} f'(x_+) - \lim_{x_- \nearrow x} f'(x_-).
\end{align*}
We have $\left[ f' \right](\xi)\neq 0$ for $\xi\in\Theta$ and $\left[ f' \right](x)=0$ otherwise. 
It is the goal to construct an approximation of $\left[ f' \right](x)$ 
that rapidly converges to zero away from $\xi\in\Theta$. Consider a point $x\in I$ and a stencil $\Stencil = \{x_0,...,x_{m+2}\}$ 
of its nearest $m+3$ grid points with $x_0 < ... < x_{m+2}$, $m>1$, and known values $f(x_i)$. We define
\begin{align*}
\Stencil^+ :=&\; \left\{ x_i : x_i\ge x\right\}\; \text{ and } \; \Stencil^- := S_x \setminus \Stencil^+.
\end{align*}
and assume $\left|\Stencil^+\right|,\left|\Stencil^-\right|>1$ to have sufficient information left and right of $x\in I$. Let
$k$ with $1\le k<m+1$ such that $x_j \in \Stencil^-$ for $j\le k$ and $x_j \in \Stencil^+$ for $j>k$. Further, define the
maximum separation length $h_x:=\max\{|x_i-x_{i-1}|:x_i,x_{i-1}\in S_x\}$. An approximation of $\left[ f' \right](x)$ is
now defined by
\begin{align*}
L_{m,\Stencil} f(x) = h_x^{m-1} \sum_{i=0}^{m+2} c_i f(x_i),
\end{align*}
where the coefficients $c_i$ are uniquely determined through the system of linear equations,
see \cite[Theorem 2.1]{ArchibaldGelbYoon2008},
\begin{align*}
 c_0 p_l(x_0) + \cdots + c_{m+2} p_l(x_{m+2})  =&\; p_l^{(m)}(x),\quad l=0,...,m,\\[2mm]
 c_{k+1} p_l(x_{k+1}) + \cdots + c_{m+2} p_l(x_{m+2}) =&\; h_x^{1-m} p_l^{(1)}(x),\quad l=0,1.
\end{align*}
Here, $p_0,...,p_m$ is a basis of the space of all polynomials up to degree $m$ and
$p_l^{(k)}$ denotes the $k$-th derivative of $p_l$. We choose $p_l(x)=x^l\in\Pi_m$ as 
basis and apply the algorithm from \cite[Appendix A]{ArchibaldGelbYoon2008}
to calculate the coefficients $c_i$, $i=0,\ldots,m+2$. Note that $c_i=\gO(h_x^{-m})$. 

In \cite[Corollary 2.5]{ArchibaldGelbYoon2008}, it was shown that for the special case when
$x_i\le x,\xi\le x_{i+1}$ for some $x_i\in \Stencil$, i.e., the discontinuity point $\xi$ and 
the reconstruction point $x$ lie within the same cell, and $\xi$ being the only jump discontinuity in the
smallest closed interval $I_x$ that contains $S_x$, then it holds
\begin{align*}
L_{m, \Stencil} f(x) =&\; \left[ f' \right](\xi) + \gO(h_x^{ \min(\kappa,m)})
\end{align*}
for piecewise smooth $f\in C^\kappa(I_x\setminus \{\xi\})$, $\kappa>1$. Hence, the approximation $L_{m, \Stencil} f(x)$
converges to $\left[ f' \right](\xi)$  with a rate depending on $m$ and the local smoothness of $f$.
\begin{remark}
As explained in \cite{ArchibaldGelbYoon2008} in a heuristic manner,
the derivative discontinuity detection method is based on a polynomial annihilation technique. If
the function $f$ is smooth in $I_x$, then $L_{m, \Stencil} f(x)$ will annihilate all terms in
the Taylor expansion of $f$ up to degree $\kappa -1$ by construction of the $c_i$. Since the residual of 
the Taylor expansion is bounded, multiplication with the scaling factor $h_x^{m-1}$ produces
an error which is proportional to the density of the sampled points. Hence, there is a high order
of convergence in smooth regions, whereas function values on both sides of $x$, i.e., in
$\Stencil^-$ and $\Stencil^+$, allow the approximation of the derivative discontinuity in the
point $x$ when $f$ has a kink there.
\end{remark}
In the setting of knot trees, a slight modification is necessary at the boundary. Let $x_0<x$ 
be the only point available in $\Stencil^-$ and $p$ the interpolating polynomial in $\Stencil[x]$ with
$|\Stencil|=3$, i.e., two further knots right of $x$ are available, then we define 
\begin{align*}
L_{\Stencil} f(x) :=&\; p'(x) - \frac{p(x)-f(x_0)}{x-x_0}
\end{align*}
as approximation of $\left[ f' \right](x)$. The analogous expression is used at the right 
boundary.

We will now explain the use of the kink detection in our $hp$-GSG method. It is important
to identify a kink as fast as possible in order to select suitable basis functions. Hence,
we apply the approach in the subroutine \textit{Selection} of Algorithm~\ref{AlgoSparseGridMain}.
Suppose $\knotHier[\bm{l}][\bm{j}]$ was added by changing $\knotHier[l'_d][j'_d]$ in its parent 
knot $\knotHier[\bm{l'}][\bm{j'}]$ along the $d$-axis.
For levels $q\le 2$, no kink detection is performed since there are only a few points 
available. If $q>2$, we compute $\eta_d:=|L_{m, \Stencil} f(\knotHier[l_d][j_d])|$ using the $m+2:=4$
nearest knots to $x:=\knotHier[l_d][j_d]$ along the $d$-axis and the knot itself such that
$x_0<x_1<x<x_3<x_4$. If $\eta_d>w_{kink}$
with a certain threshold $w_{kink}$, we conclude that a kink is detected in the one-dimensional support 
$\support[l_d][j_d]$ and hence we choose the piecewise linear basis function $\baseFuncHier[l_d][j_d,1]$
there. Otherwise, we set $\support[l_d][j_d]=\baseFuncHier[l_d][j_d,p_d]$ with 
$p_d:=\min\{p'_d+1,p_{\max}\}$, where $p'_d$ is the degree of the basis function of the corresponding
parent knot. It might happen in our adaptive algorithm that there are not enough points along a
certain axis to perform a kink detection. However, in this case we expect the interpolation error to be 
very small there and again choose the piecewise linear basis function.

Eventually, the subroutine \textit{Modification} is skipped, i.e., no modification for 
parent knots takes place.  For later use, we call this method $hp$-GSG-k, where \glqq k\grqq{} stands 
for \textbf{k}ink detection.

\section{Numerical examples}\label{sec:num}
We will investigate the performance of the newly proposed $hp$-GSG methods compared to three already
existing approaches: choosing always (i) locally linear basis functions as in \cite{MaZabaras2009}, 
(ii) locally quadratic basis functions with anisotropic inspection of the sparse
grid indices as favoured in 
\cite{JakemanRoberts2011} for non-smooth functions, and
(iii) highest degree basis functions first proposed in \cite{Bungartz1998}, referred 
to as \glqq linear\grqq{}, \glqq quadratic\grqq{} and \glqq highest\grqq{}. 
The strategies \glqq linear\grqq{} and \glqq highest\grqq{} are special cases
of Algorithm~\ref{AlgoSparseGridMain}, whereas programming \glqq quadratic\grqq{} requires some adaptation. 
In the latter case, the threshold $w_{max}$ refers to the $\varepsilon$ that controls the
error indicators in \cite[Formula (8)]{JakemanRoberts2011}. The algorithms have been
programmed using C++ and Matlab. We have used the Matlab-Version 
R2024a on an Intel Xeon Gold 6130 CPU. For time measurements, an average over three runs was taken.

The following functions which have been considered
in previous work \cite{JakemanRoberts2011,MaZabaras2009,TaoJiangCheng2021} are used:
\begin{align*}
\fcurve[] (x) =&\; \frac{1}{\left| 0.3 - x_1^2 - x_2^2 \right| + 0.1}, \; x \in [0,1]^2,\\[2mm]
\fc (x) =&\; \exp\left( - \sum_{i=1}^{N} a_i \left| x_i - 0.51 \right| \right), \; x \in [0,1]^N,
\;a_i=2^{3-i},i=1,\ldots,N,\\[2mm]
\fg (x) =&\; \prod_{i=1}^N \frac{4\,| x_i^2- 0.66^2 | + a_i }{a_i + 1}, \; x \in [0,1]^N,
\;a_1=0.5,\,a_i=(i-1)^2,i=2,\ldots,N,\\[2mm]
\fd (x)=&\; 
\begin{cases}
0 &\text{if } \max\{x_1, x_2\} > 0.51 \\
\exp\left( \sum_{i=1}^{N} a_i x_i \right) & \text{else}
\end{cases}, \; x \in [0,1]^N, a_i=2^{3-i},i=1,\ldots,N.
\end{align*}
The first problem is a well established benchmark problem with a one-dimensional kink that is
not aligned with a certain coordinate axis. Clearly, the standard dimension-adaptive 
(anisotropic) sparse grid method fails in this case since it cannot resolve such non-smooth 
functions, see the discussion in \cite[Chapter 4.1]{MaZabaras2009}. The function $\fc (x)$ is taken from 
the continuous family of the Genz library \cite{Genz1984}. We set for the dimension $N=2,10,50$.
Eventually, we consider a modified g-function $\fg (x)$ of Sobol \cite{SaltelliSobol1995} with
quadratic terms $x_i^2$ instead of linear ones and set $N=2,5$. As a robustness check,
we also consider a discontinuous function $\fd (x)$ from the Genz library \cite{Genz1984} for $N=2,5$. 
In all cases, we choose the parameters $a_i$ in such a way to mimic a decreasing 
importance of the directions $x_i$. Note that we use $0.51$ in the definition of 
$\fc$ and $\fd$ for the characterization of the kinks and discontinuities, respectively, to avoid
a direct alignment with the knots.

The error is measured in the following norms:
\begin{align*}
\errMax &:= \max_{i=1,...,M} \left| f(x^{(i)}) - \mathscr{U}_{A}^{(q)}(f) (x^{(i)}) \right|, \\[2mm]
\errTwo &:= \left( \frac{1}{M}\sum_{i=1}^M \left| f(x^{(i)}) - \mathscr{U}_{A}^{(q)}(f) (x^{(i)}) 
\right|^2 \right)^{1/2},
\end{align*}
with randomly chosen test points $x^{(i)}\in\Omega$, $i=1,...,10^5$, and
additional $10^3$ points $x^{(j)}$ along the kinks. Here, $f$ and $\mathscr{U}_{A}^{(q)}(f)$
denote the true function and its interpolant, respectively. 

We set $p_{\max}=6$, $q_{\min}=1$, and $q_{\max}=25$ in Algorithm~\ref{AlgoSparseGridMain} for 
all numerical tests, except for $\fcurve[]$, where we use $q_{\max}=30$.
\begin{figure}[ht!]
\centering
\includegraphics[width=7cm]{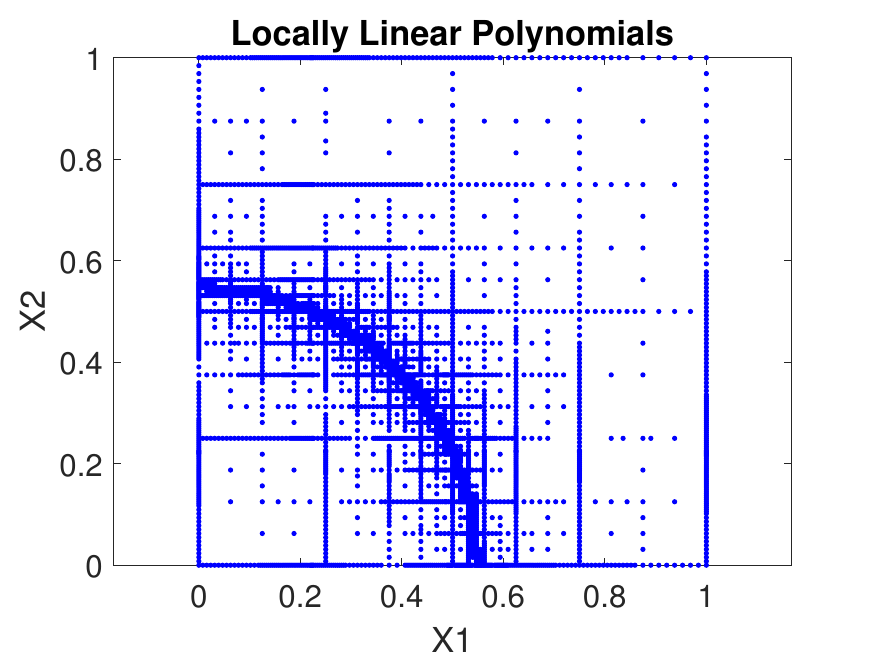}
\hspace{0.1cm}
\includegraphics[width=7cm]{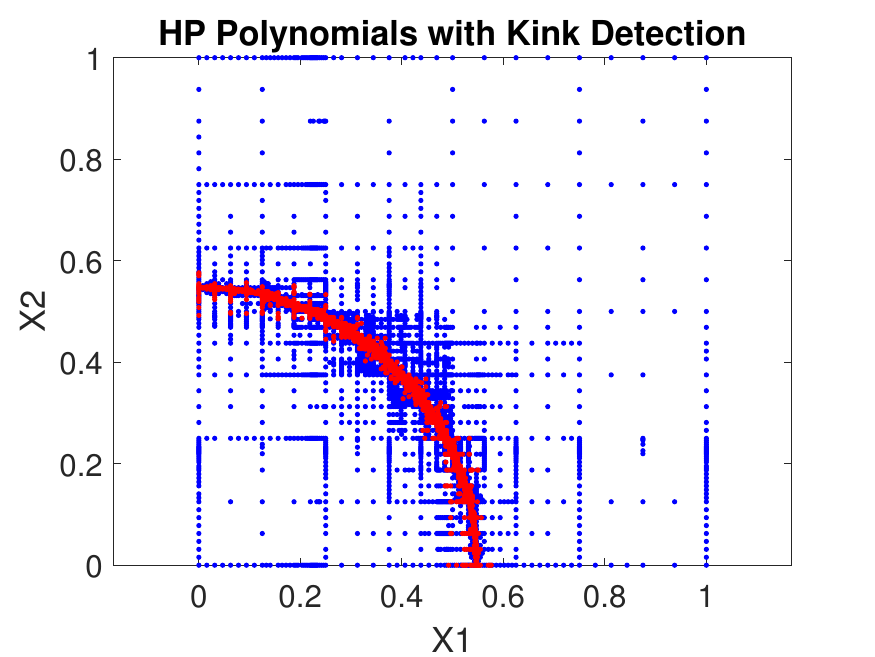}
\parbox{14cm}{
\caption{Test problem with $\fcurve[]$. Knot placement for locally linear 
polynomials (left) and $hp$-GSG-k (right) with thresholds $w_{\max}=10^{-3}$
and $w_{kink}=5$. Knots of basis functions for which a kink was detected 
within an interval of size smaller than $2^{-6}$ are colored in red.}
\label{fig:2d-points}
}
\end{figure}
\begin{figure}[ht!]
\centering
\includegraphics[width=7cm]{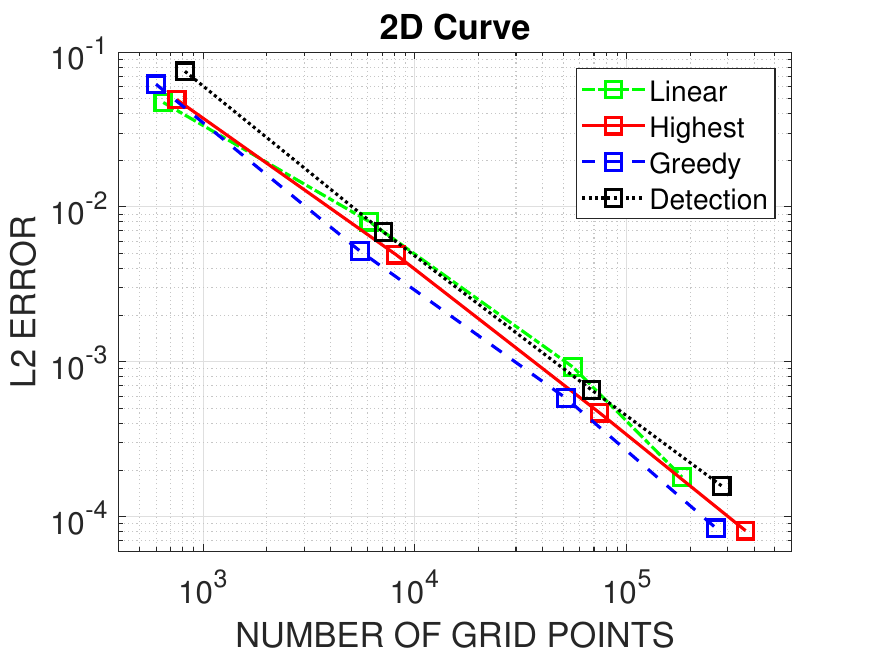}
\hspace{0.1cm}
\includegraphics[width=7cm]{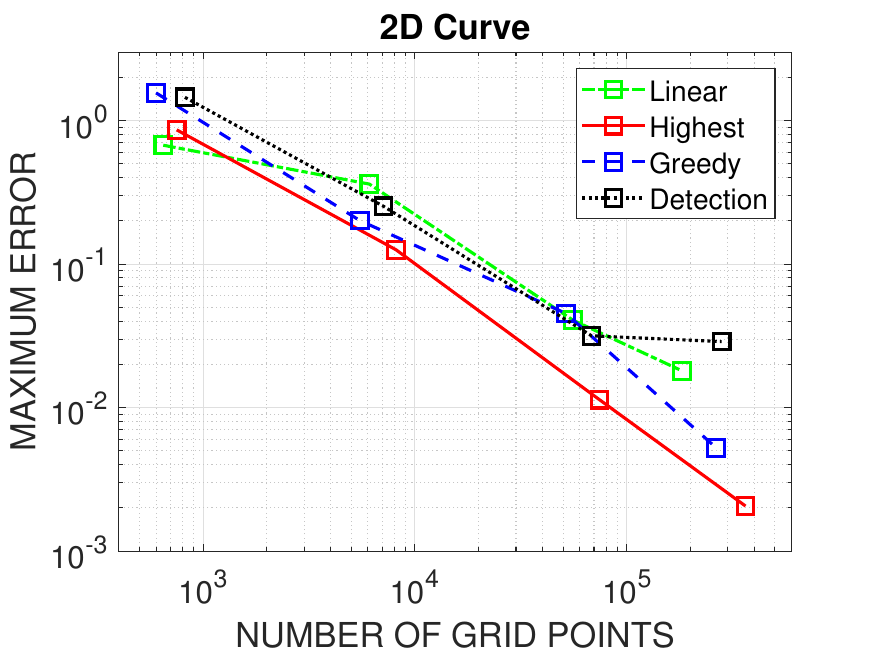}
\parbox{14cm}{
\caption{Test problem with $\fcurve[]$. Comparison of the errors $\errTwo$ (left) and
$\errMax$ (right).}
\label{fig:2d-errors}
}
\end{figure}
\subsection{Approximation of $\fcurve[]$}
The function $\fcurve[]$ exhibits a regularized line singularity along the circle around $(0,0)$ with
radius $r=\sqrt{0.3}$. The gradients close to this line are quite large, which motivates a larger
threshold $w_{kink}=5$ for the kink detection. We use 
$w_{max}=10^{-2i+1}$, $i=2,\ldots,6$, for \glqq quadratic\grqq{}  
and and $w_{\max}=10^{-i}$, $i=1,\ldots,4$,
for all other methods. In Figure~\ref{fig:2d-points}, the knot placement for locally linear 
polynomials (method \glqq linear\grqq{}) and $hp$-GSG-k with thresholds $w_{\max}=10^{-3}$ is shown. The
kink detection is able to identify the region of interest and places less grid points in smooth regions
due to the higher-order polynomials used there. However, the advantage of the $hp$-strategy in terms of
accuracy is less pronounced as can be seen from the error plots in Figure~\ref{fig:2d-errors}. All methods
perform quite similar for the 2-norm, whereas the simple strategy \glqq highest\grqq{} and the
more sophisticated \glqq quadratic\grqq{} are best 
for the $\infty$-norm closely followed by $hp$-GSG-g. This example verifies that methods with
local refinement deliver comparably accurate results. A closer look at the computing time shows that all
methods perform equally well with respect to achieved accuracy. Therefore, we omit the details.
\begin{figure}[ht!]
\centering
\includegraphics[width=7cm]{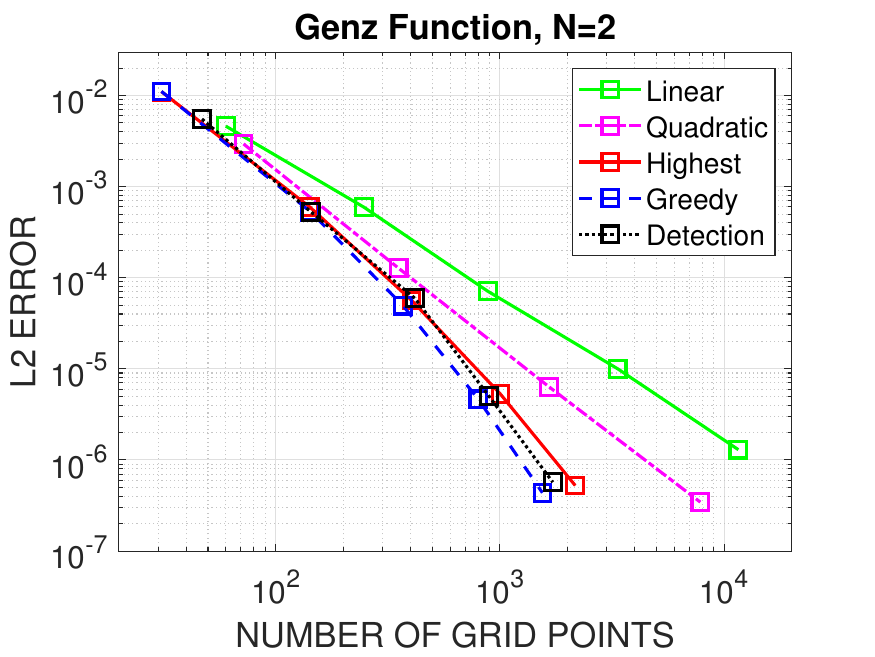}
\hspace{0.1cm}
\includegraphics[width=7cm]{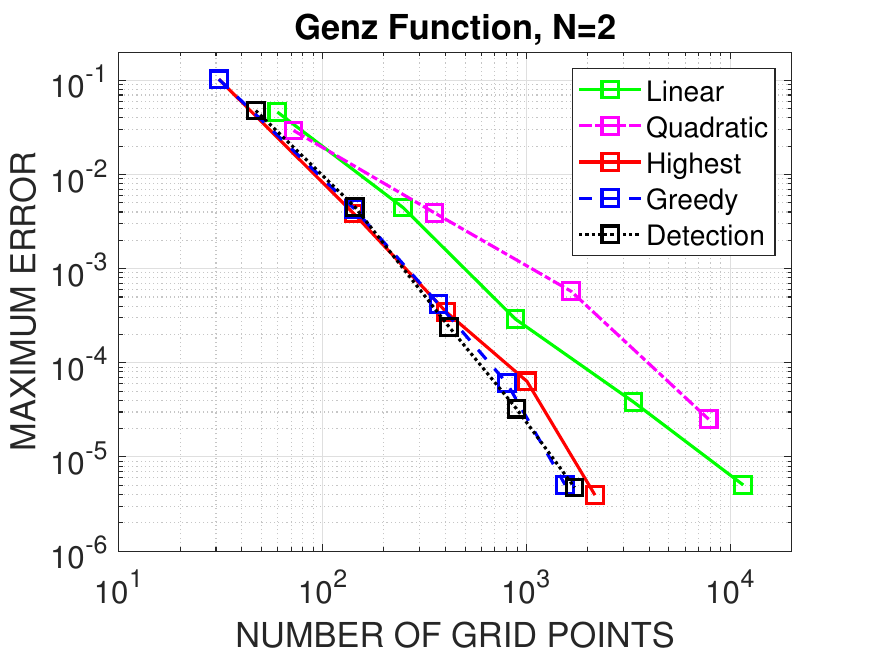}
\parbox{14cm}{
\caption{Test problem with $\fc[2]$. Comparison of the errors $\errTwo$ (left) and
$\errMax$ (right).}
\label{fig:c0func-d2-errors}
}
\end{figure}
\begin{figure}[ht!]
\centering
\includegraphics[width=7cm]{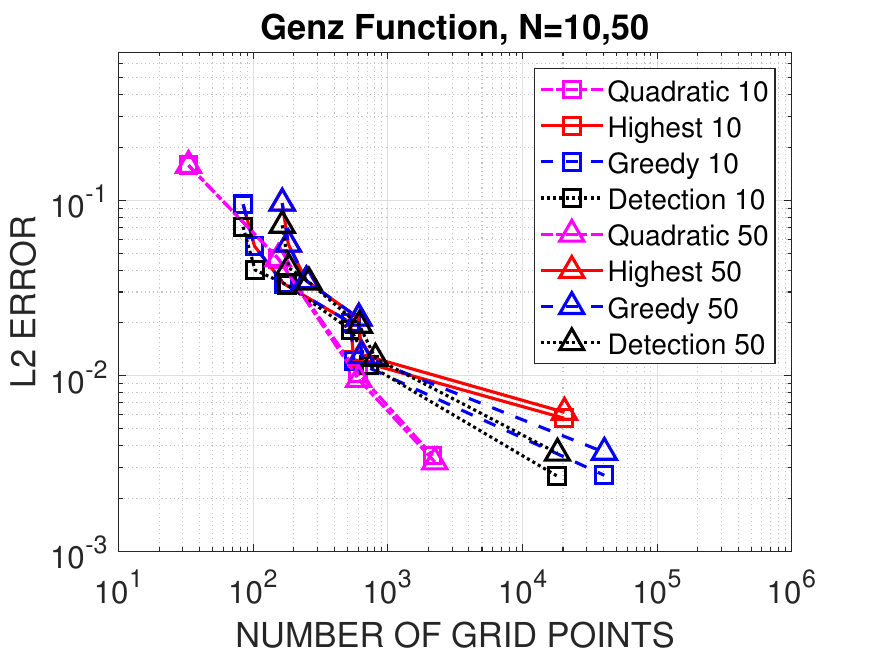}
\hspace{0.1cm}
\includegraphics[width=7cm]{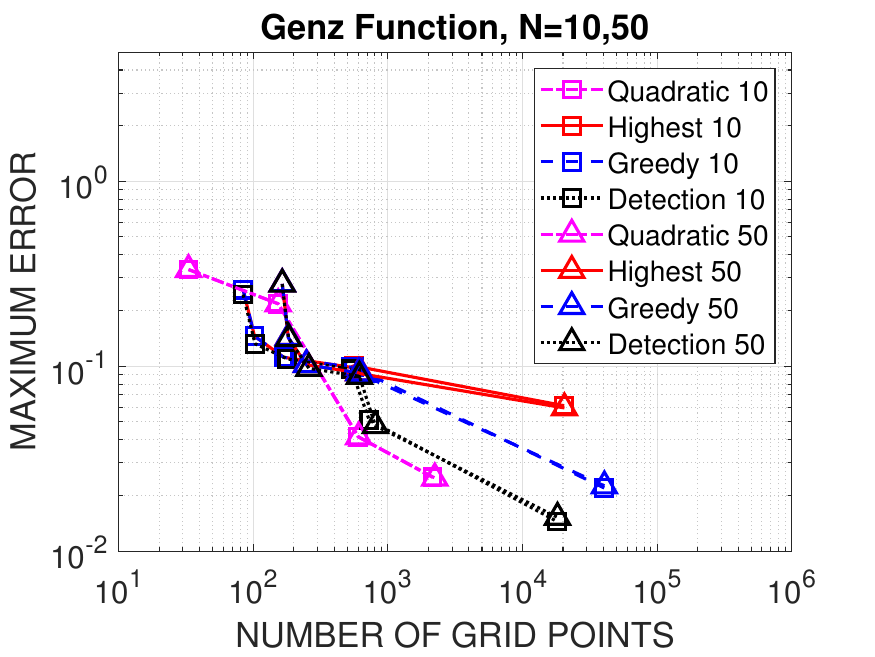}
\parbox{14cm}{
\caption{Test problem with $\fc[10,50]$. Comparison of the errors $\errTwo$ (left) and
$\errMax$ (right).}
\label{fig:c0func-d10a50-errors}
}
\end{figure}
\subsection{Approximation of $\fc$ in various dimensions}
\begin{figure}[ht!]
\centering
\includegraphics[width=7cm]{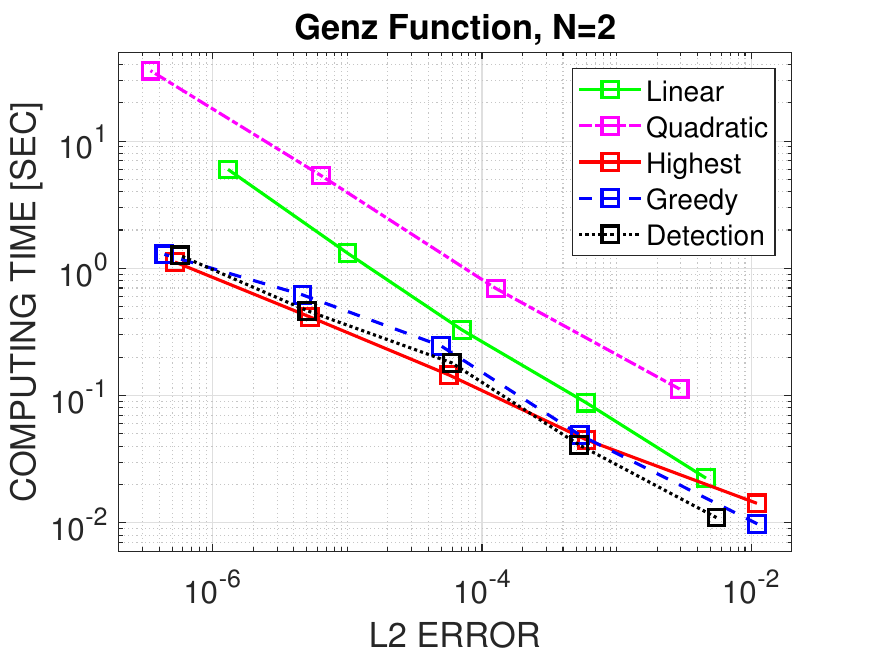}
\hspace{0.1cm}
\includegraphics[width=7cm]{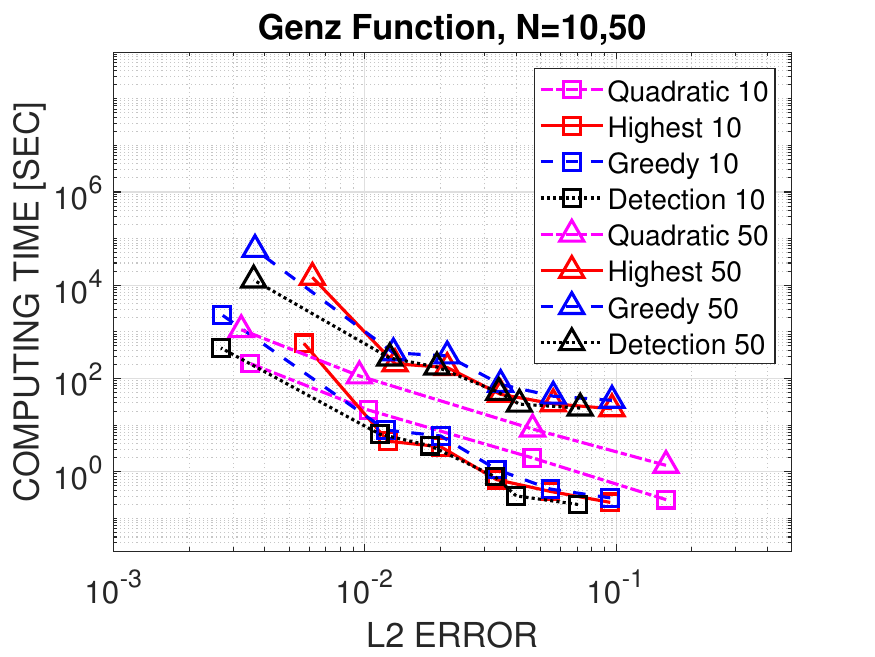}
\parbox{14cm}{
\caption{Test problem with $\fc[N]$. Comparison of computing time in seconds versus accuracy in
the 2-norm for $N=2$ (left) and $N=10,50$ (right).}
\label{fig:c0func-cpu}
}
\end{figure}

Now we consider the function $\fc$ which has kinks of decreasing importance along the 
coordinate axis. For dimension $N=2$, we use 
$w_{max}=10^{-4},\,6\cdot 10^{-7},\,3\cdot 10^{-9},\,10^{-11},$ for \glqq quadratic\grqq{}  
and a sequence of $5$ thresholds $w_{\max}=10^{-i}$, $i=2,\ldots,6$, for all 
other methods. For $N=10,50$, we choose
$w_{\max}=10^{-i}$, $i=2,\ldots,5$, for \glqq quadratic\grqq{} and 
$w_{\max}\in [3\cdot 10^{-3},10^{-1}]$ for $6$ different runs else. The threshold for 
the kink detection is fixed
with $w_{kink}=1$. The corresponding convergence results are
plotted in Figure~\ref{fig:c0func-d2-errors} and Figure~\ref{fig:c0func-d10a50-errors}.
For $N=2$, both newly designed $hp$-GSG methods perform better than the other three
for both error measures, $\errTwo$ and $\errMax$. Further, the additional effort
of \glqq quadratic\grqq{} to detect anisotropic features does not pay off.
The strategy \glqq linear\grqq{} is no longer able to compete, which becomes even more 
evident for higher dimensions. We
therefore omit its results for $N=10,50$. In these cases, the new methods outperform
\glqq highest\grqq{} for higher accuracies. For higher tolerances, 
\glqq quadratic\grqq{} is the best method. However, in terms of computing time
versus the accuracy achieved in the 2-norm, it is less efficient for $N=2,10$, 
as can be seen in Figure~\ref{fig:c0func-cpu}. The new approaches work equally well for $N=2$.
They still outperform \glqq highest\grqq{} in higher dimensions, but 
they are less efficient for higher accuracies and $N=50$ than \glqq quadratic\grqq{}.
The ability to detect anisotropies is very useful here. Similar observations are made for 
the $\infty$-norm.
\subsection{Approximation of $\fg$ in various dimensions}
\begin{figure}[ht!]
\centering
\includegraphics[width=7cm]{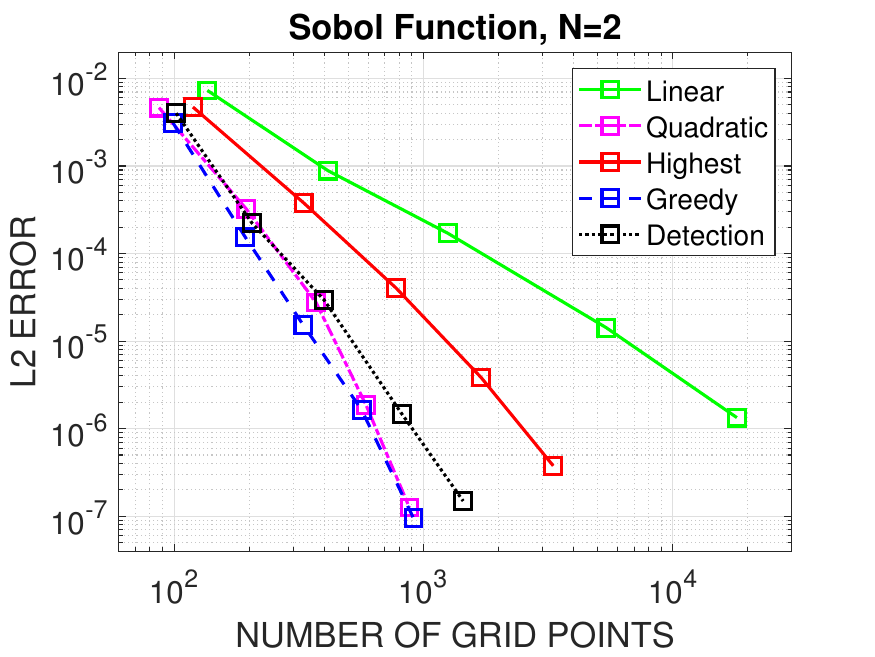}
\hspace{0.1cm}
\includegraphics[width=7cm]{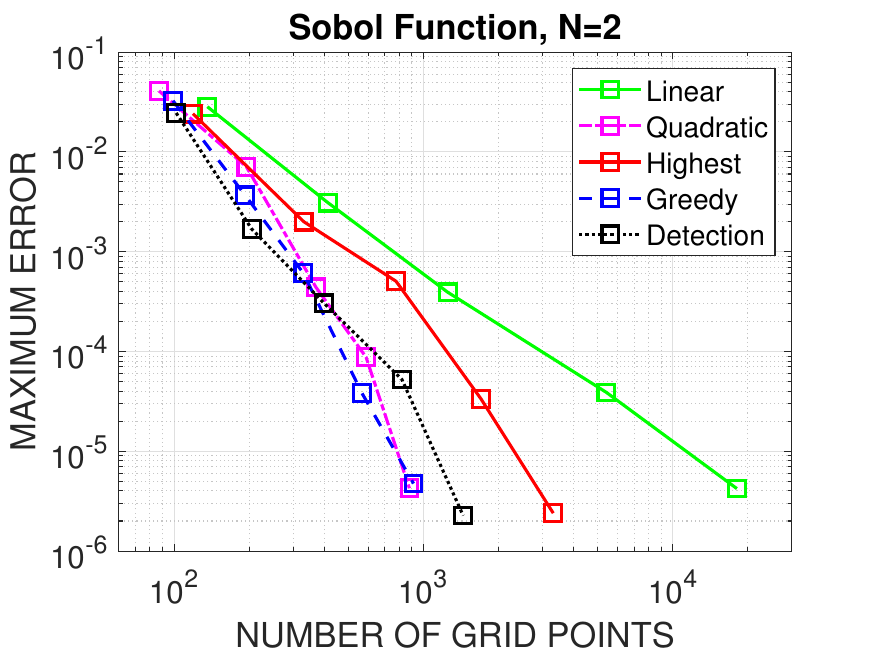}
\parbox{14cm}{
\caption{Test problem with $\fg[2]$. Comparison of the errors $\errTwo$ (left) and
$\errMax$ (right).}
\label{fig:gpfunc-d2-errors}
}
\end{figure}

\begin{figure}[ht!]
\centering
\includegraphics[width=7cm]{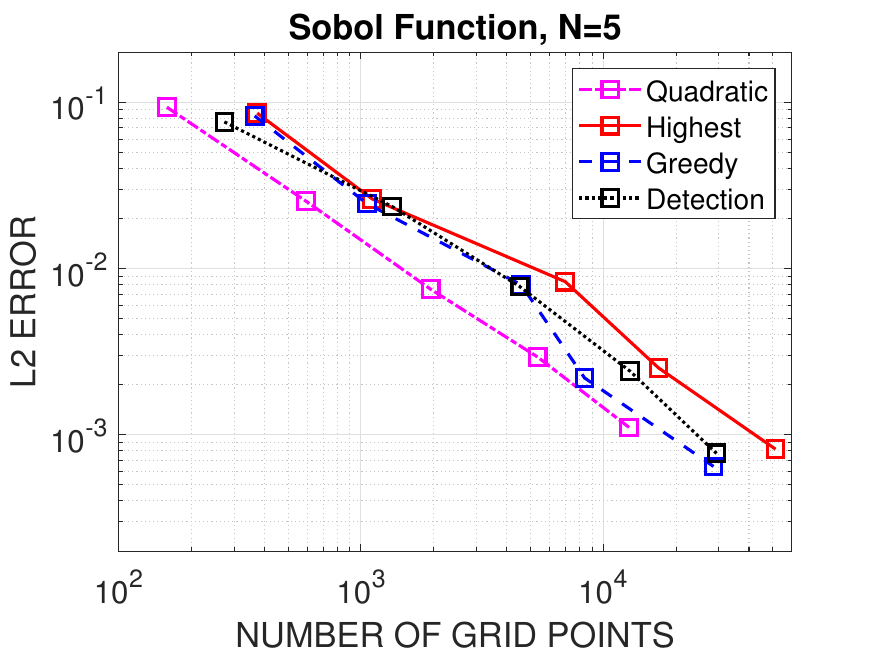}
\hspace{0.1cm}
\includegraphics[width=7cm]{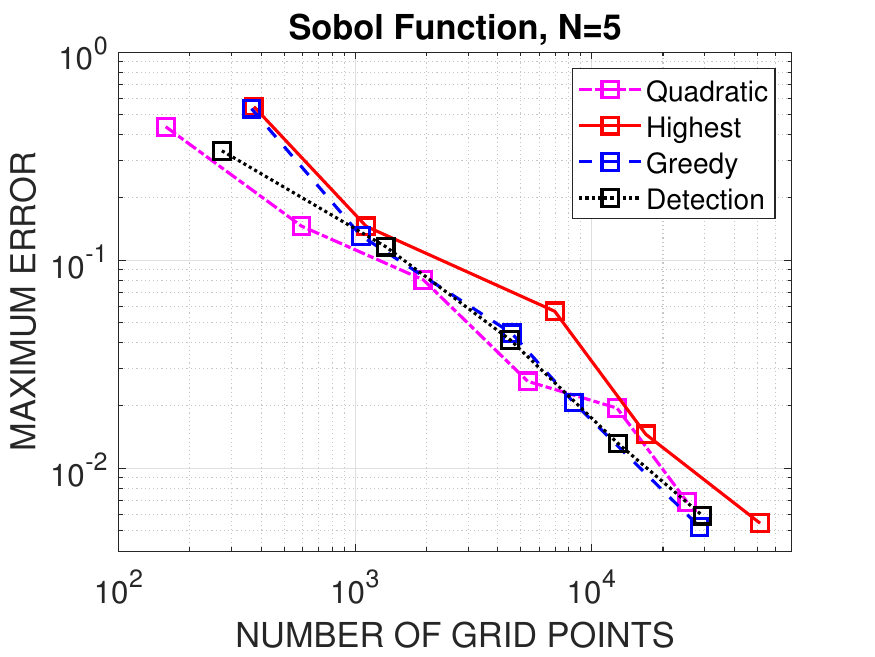}
\parbox{14cm}{
\caption{Test problem with $\fg[5]$. Comparison of the errors $\errTwo$ (left) and
$\errMax$ (right).}
\label{fig:gpfunc-d5-errors}
}
\end{figure}

\begin{figure}[ht!]
\centering
\includegraphics[width=7cm]{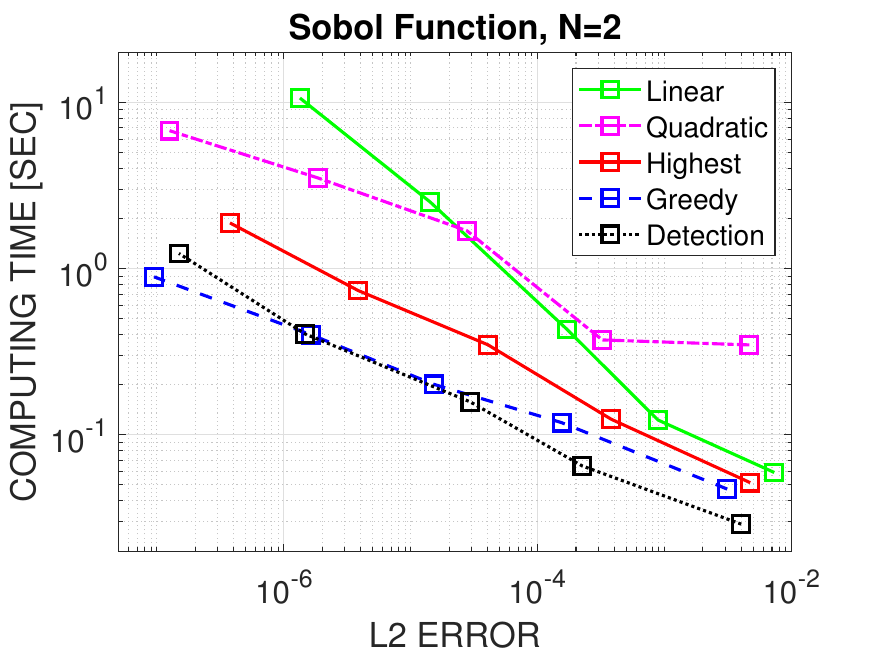}
\hspace{0.1cm}
\includegraphics[width=7cm]{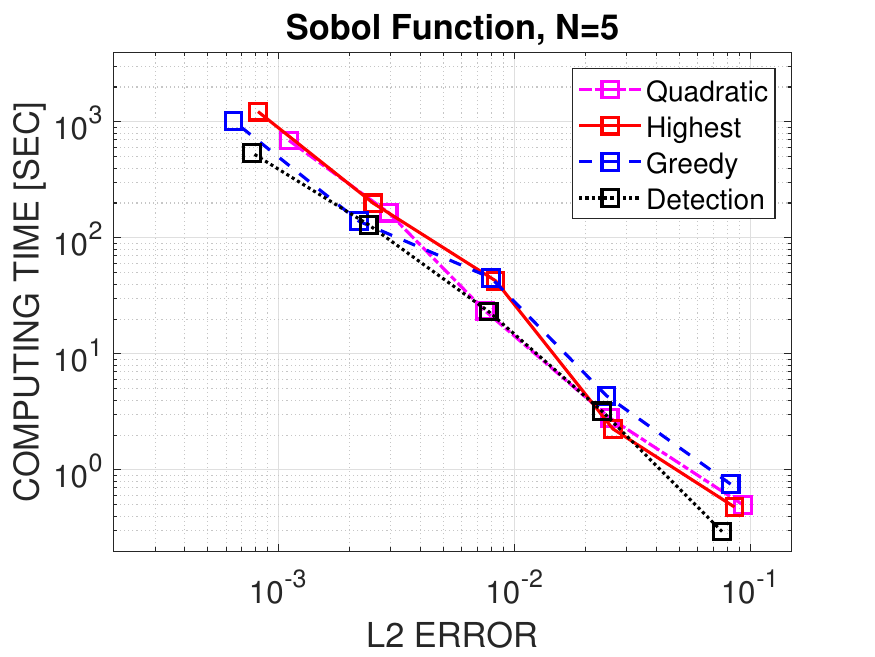}
\parbox{14cm}{
\caption{Test problem with $\fg[N]$. Comparison of computing time in seconds versus accuracy in
the 2-norm for $N=2$ (left) and $N=5$ (right).}
\label{fig:gpfunc-cpu}
}
\end{figure}

Next we consider the function $\fg$ for $N=2,5$. For each index $i$, a lower value of 
$a_i$ indicates a higher importance of the variable $x_i$. The threshold $w_{\max}$
for \glqq quadratic\grqq{} varies in $[10^{-13},10^{-4}]$ for $N=2$ and is set to
$10^{-i}$, $i=3,\ldots,7,$ for $N=5$. For the other methods, we use 
$w_{max}=10^{-i}$, $i=2,\ldots,6$, for $N=2$ and five varying values in $[10^{-3},10^{-1}]$ for $N=5$.
The threshold for the kink detection is again fixed with $w_{kink}=1$. In 
Figure~\ref{fig:gpfunc-d2-errors} and Figure~\ref{fig:gpfunc-d5-errors}, the
convergence results are shown for $\errTwo$ and $\errMax$. For $N=2$, the greedy approach
$hp$-GSG-g and \glqq quadratic\grqq{} perform best, whereas for $N=5$, $hp$-GSG-k 
is equally good. The advantage
of these methods compared to \glqq linear\grqq{} and \glqq highest\grqq{} is clearly
visible and quite impressive for $N=2$ and higher accuracies. Once again, the strategy 
\glqq linear\grqq{} needs too much knots for $N=5$. Therefore, the corresponding results
are not shown in Figure~\ref{fig:gpfunc-d5-errors}. The superiority of the new methods for $N=2$
is clearly reflected in the computing time shown in Figure~\ref{fig:gpfunc-cpu} as function
of the $2$-error. For $N=5$, all methods except \glqq linear\grqq{} perform quite 
the same. The same holds true for the maximum error.

\subsection{Approximation of $\fd$ in various dimensions}
\begin{figure}[t!]
\centering
\includegraphics[width=7cm]{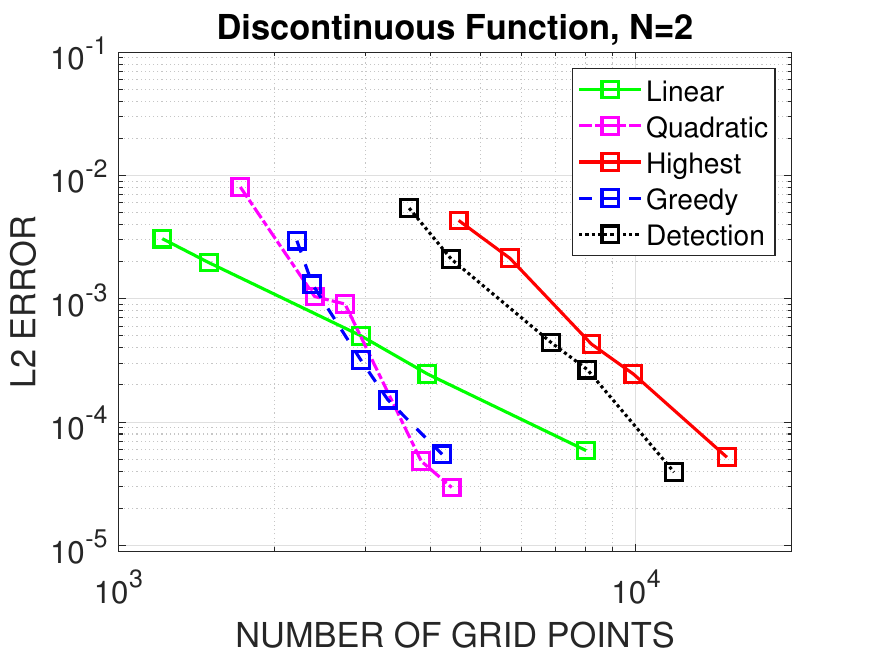}
\hspace{0.1cm}
\includegraphics[width=7cm]{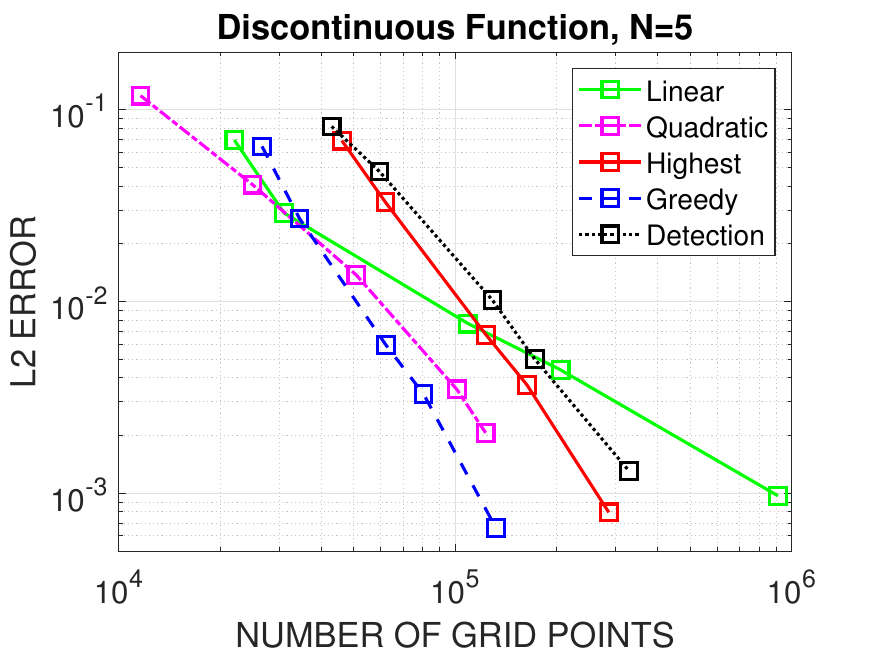}
\parbox{14cm}{
\caption{Test problem with $\fd$. Comparison of the errors $\errTwo$ 
for $N=2$ (left) and $N=5$ (right).}
\label{fig:dcfunc-l2err}
}
\end{figure}

\begin{figure}[t!]
\centering
\includegraphics[width=7cm]{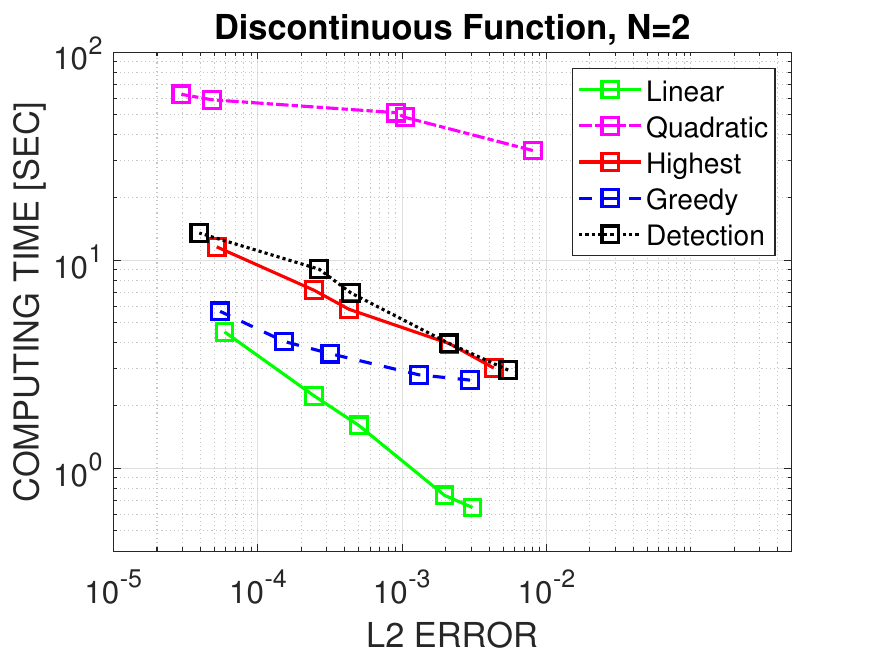}
\hspace{0.1cm}
\includegraphics[width=7cm]{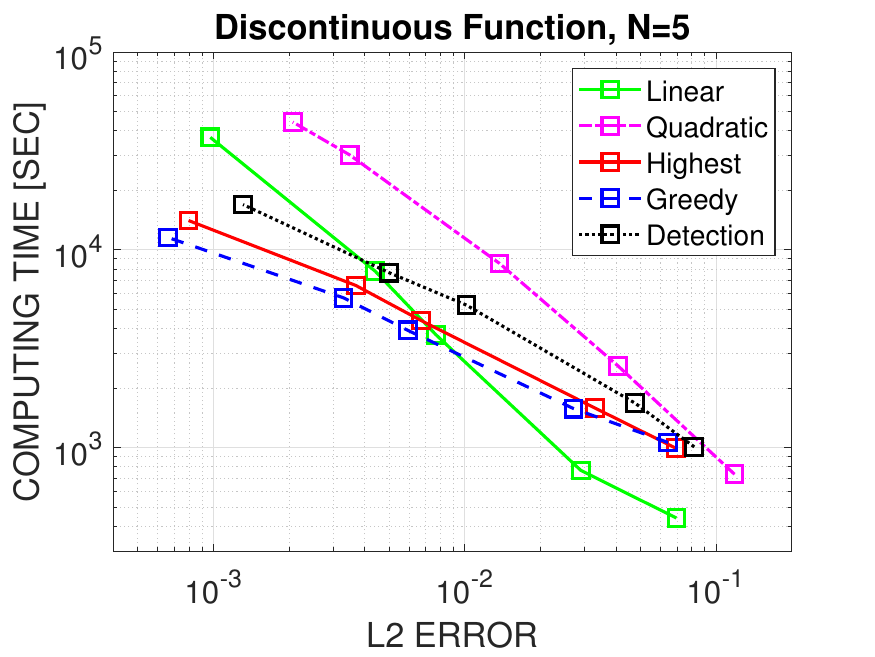}
\parbox{14cm}{
\caption{Test problem with $\fd$. Comparison of computing time versus accuracy in
the 2-norm for $N=2$ (left) and $N=5$ (right).}
\label{fig:dcfunc-cpu}
}
\end{figure}

As a test of robustness and application to discontinuous functions, we finally consider
$\fd[N]$ for $N=2,5$. The function has discontinuities along the lower-dimensional manifolds
$x_1=0.51$ and $x_2=0.51$. Its approximation by continuous polynomials is characterized by the
so-called Gibbs phenomenon which manifests itself in undesired over- and undershoots in the
neighborhood of the discontinuities. Adaptivity
and the use of linear polynomials close to the discontinuity can mitigate this effect. Although
not designed for this case, we apply our kink detection approach with $w_{kink}=1$ and without
any further adaptation. The coefficients $a_i$ are monotonically decreasing, which renders higher 
dimensions less important. The threshold $w_{\max}$
for \glqq quadratic\grqq{} varies in $[5\cdot 10^{-9},5\cdot 10^{-7}]$ and $[5\cdot 10^{-9},10^{-5}]$ 
for $N=2$ and $N=5$, respectively. For the other methods, we use the thresholds
$w_{max}=10^{-2}, 5\cdot 10^{-3}, 10^{-3}, 5\cdot 10^{-4}, 10^{-4},$ for $N=2$ and these values
multiplied by $10$ for $N=5$. In Figure~\ref{fig:dcfunc-l2err} and Figure~\ref{fig:dcfunc-cpu}, we
present results for the $2$-norm only, since the $\infty$-norm is not appropriate here. The simple
method \glqq linear\grqq{} performs quite well. It even shows quadratic order with respect to
the reciprocal number of grid points in the case $N=2$. For $N=5$, it drops down to first order. 
However, in both test cases, the method is quite efficient with respect to computing time. 
Although not designed for discontinuous functions, $hp$-GSG-k delivers acceptable results and
can keep up with \glqq highest\grqq{}. Both the greedy approach and \glqq quadratic\grqq{} perform 
best and equally well with a small but visible advantage for the first one for $N=5$. 
However, \glqq quadratic\grqq{} is clearly outperformed in terms of computing 
times by all other methods. The reason is as follows: Its exploration of indices and identification 
of admissible points become more and more expensive for higher refinement levels which are required
to accurately resolve the discontinuities.


\section{Conclusion}\label{sec:concl}
In this paper, we have developed a novel $hp$-adaptive sparse grid collocation
algorithm with spatial refinement for piecewise smooth functions with kinks along certain manifolds.
Both the greedy and the kink detection approach perform equally well for three
benchmark problems with continuous functions. In these cases, classical (isotropic) 
$h$-methods with multi-linear hierarchical basis functions 
or fixed maximum polynomial degree are clearly outperformed for higher dimensions with
respect to the number of grid points necessary to achieve a certain accuracy as well as
the corresponding computing time. The comparatively simple greedy strategy to select 
optimal local polynomial degrees works very efficient and robust. The kink detection
algorithm has the potential for further improvements as can be seen for the high-dimensional
Genz problem, but comes with the burden of an additional appropriate choice of the threshold 
parameter $w_{kink}$. The anisotropic $h$-method with quadratic basis functions 
performs well in terms of accuracy per number of grid points, but often needs significantly
higher computing time to detect important components and to check for admissible points. Its 
advantage in high dimensions with
strong anisotropy is obvious as well as its deficiencies for the discontinuous function considered. Here,
the $hp$-greedy method combines the merit of using linear polynomials in the neighborhood of
the discontinuity and higher order polynomials elsewhere.  

In future work, we will apply our stochastic collocation method to real gas network problems,
where regulations of the gas flow yields kinks in certain quantities of interest as e.g.
the maximum pressure over time at delivery exits. Efficient approximations of such quantities
in higher dimensions will be one of the key steps to solve stochastic optimal control problems with
probabilistic constraints.  

\vspace{0.5cm}
\par
\noindent {\bf Acknowledgements.}
Both authors are supported by the Deutsche Forschungsgemeinschaft
(German Research Foundation) within the collaborative research center
TRR154 {\em ``Mathematical modeling, simulation and optimisation using
the example of gas networks''} (Project-ID 239904186, TRR154/3-2022, TP B01).\\

\par
\noindent {\bf Declaration of competing interest.}
The authors declare that they have no known competing financial interests or personal 
relationships that could have appeared to influence the work reported in this paper.

%------------------
\bibliographystyle{plain}
\bibliography{bibhpcoll}

\begin{thebibliography}{10}

\bibitem{ArchibaldGelbYoon2005}
R.~Archibald, A.~Gelb, and J.~Yoon.
\newblock Polynomial fitting for edge detection in irregularly sampled signals
  and images.
\newblock {\em SIAM J. Numer. Anal.}, 43:259--279, 2005.

\bibitem{ArchibaldGelbYoon2008}
R.~Archibald, A.~Gelb, and J.~Yoon.
\newblock Determining the locations of discontinuities in the derivatives of
  functions.
\newblock {\em Appl. Numer. Math.}, 58:577--592, 2008.

\bibitem{Bungartz1998}
H.-J. Bungartz.
\newblock {\em Finite elements of higher order on sparse grids}.
\newblock Habilitationsschrift. Institut f\"ur Informatik, TU M\"unchen and
  Shaker Verlag, Aachen, 1998.

\bibitem{BungartzGriebel2004}
H.-J. Bungartz and M.~Griebel.
\newblock Sparse grids.
\newblock {\em Acta Numer.}, 13:147--269, 2004.

\bibitem{FuchsGarcke2020}
B.~Fuchs and J.~Garcke.
\newblock Simplex stochastic collocation for piecewise smooth functions with
  kinks.
\newblock {\em International Journal for Uncertainty Quantification}, 10:1--24,
  2020.

\bibitem{Genz1984}
A.~Genz.
\newblock Testing multidimensional integration routines.
\newblock In B.~Ford, J.C. Rault, and F.~Thomasset, editors, {\em International
  conference on tools, methods and languages for scientific and engineering
  computation}, pages 81--94. Elsevier North-Holland, Inc., 1984.

\bibitem{GerstnerGriebel1998}
T.~Gerstner and M.~Griebel.
\newblock Numerical integration using sparse grids.
\newblock {\em Numer. Algorithms}, 18:209--232, 1998.

\bibitem{GerstnerGriebel2003}
T.~Gerstner and M.~Griebel.
\newblock Dimension-adaptive tensor-product quadrature.
\newblock {\em Computing}, 71:65--87, 2003.

\bibitem{Griebel1998}
M.~Griebel.
\newblock Adaptive sparse grid multilevel methods for elliptic {PDE}s based on
  finite differences.
\newblock {\em Computing}, 61:151--179, 1998.

\bibitem{JakemanRoberts2011}
J.D. Jakeman and S.G. Roberts.
\newblock Local and dimension adaptive sparse grid interpolation and
  quadrature.
\newblock Technical Report https://arxiv.org/abs/1110.0010v1, 2011.

\bibitem{Klimke2006}
A.~Klimke.
\newblock {\em Uncertainty Modeling using Fuzzy Arithmetic and Sparse Grids}.
\newblock PhD thesis, Universit\"at Stuttgart, Shaker Verlag, Aachen, 2006.

\bibitem{MaZabaras2009}
X.~Ma and N.~Zabaras.
\newblock An adaptive hierarchical sparse grid collocation algorithm for the
  solution of stochastic differential equations.
\newblock {\em J. Comput. Phys.}, 228:3084--3113, 2009.

\bibitem{NobileTemponeWebster2008}
F.~Nobile, R.~Tempone, and C.~Webster.
\newblock An anisotropic sparse grid stochastic collocation method for partial
  differential equations with random input data.
\newblock {\em SIAM J. Numer. Anal.}, 46:2411--2442, 2008.

\bibitem{ObersteinerBungartz2021}
M.~Obersteiner and H.-J. Bungartz.
\newblock A generalized spatially adaptive sparse grid combination technique
  with dimension-wise refinement.
\newblock {\em SIAM J. Sci. Comput.}, 43:A2381--A2403, 2021.

\bibitem{SaltelliSobol1995}
A.~Saltelli and I.M. Sobol.
\newblock Sensitivity analysis for nonlinear mathematical models: numerical
  experience.
\newblock {\em Matem. Mod.}, 7:16--28, 1995.

\bibitem{Smolyak1963}
S.~Smolyak.
\newblock Quadrature and interpolation formulas for tensor product of certain
  classes of functions.
\newblock {\em Soviet Math. Dokl.}, 4:240--243, 1963.

\bibitem{TaoJiangCheng2021}
Z.~Tao, Y.~Jiang, and Y.~Cheng.
\newblock An adaptive high-order piecewise polynomial based sparse grid
  collocation method with applications.
\newblock {\em J. Comput. Phys.}, 433:109770, 2021.

\bibitem{Zenger1991}
C.~Zenger.
\newblock Sparse grids.
\newblock In W.~Hackbusch, editor, {\em Parallel Algorithms for Partial
  Differential Equations}, volume~31 of {\em Notes Numer. Fluid Mech}, pages
  241--251. Vieweg, Brunswick, Germany, 1991.

\end{thebibliography}

\end{document}